\documentclass[11pt]{amsart} 
\usepackage{amsthm,amsmath,amssymb,latexsym}
\usepackage{hyperref}
\usepackage{lineno}
\usepackage{enumitem}
\usepackage{chngcntr}
\usepackage{pst-all}
\usepackage{enumitem}
\counterwithin{equation}{subsection}
\begin{document}

\oddsidemargin=18pt \evensidemargin=18pt
\newtheorem{theorem}{Theorem}[section]
\newtheorem{definition}[theorem]{Definition}
\newtheorem{proposition}[theorem]{Proposition}
\newtheorem{lemma}[theorem]{Lemma}
\newtheorem{notation}[theorem]{Notation}
\newtheorem{corollary}[theorem]{Corollary}
\newtheorem{question}{Question}
\newtheorem{fact}{Fact}[section]
\newtheorem{claim}[theorem]{Claim}
\theoremstyle{definition}
\newtheorem{examples}[theorem]{Examples}
\newtheorem{remark}[theorem]{Remark}

\author{Fran\c coise {Point}$^{(\dagger)}$}
\address{Department of Mathematics (De Vinci)\\ UMons\\ 20, place du Parc 7000 Mons, Belgium}
\email{point@math.univ-paris-diderot.fr}
\title[Definable groups in topological fields with a generic derivation]{Definable groups in topological fields with a generic derivation}
\date{\today}
\def\rem{\begin{remark}\upshape}
\def\erem{\end{remark}}
\def\nota{\begin{notation}\upshape}
\def\enota{\end{notation}}
\def\defn{\begin{definition}\upshape}
\def\edefn{\end{definition}}
\def\cl{\begin{claim}\noindent\upshape}
\def\ecl{\end{claim}\par\noindent\upshape}
\def\thm{\begin{theorem}}
\def\ethm{\end{theorem}}
\def\lmm{\begin{lemma}}
\def\elmm{\end{lemma}}
\def\qed{\hfill$\quad\Box$}
\def\pr{\par\noindent{\em Proof: }}
\def\prcl{\par\noindent{\em Proof of Claim: }}
\def\kor{\begin{corollary}}
\def\ekor{\end{corollary}}
\def\frag{\begin{question}\upshape}
\def\efrag{\end{question}}
\def\prop{\begin{proposition}}
\def\eprop{\end{proposition}}
\def\fct{\begin{fact}}
\def\efct{\end{fact}}
\def\ex{\begin{example}}
\def\eex{\end{example}}

\newcommand{\QED}{{\unskip\nobreak\hfil\penalty50%
\hskip1em\hbox{}\nobreak\hfil $\Box$%
\parfillskip=0pt \finalhyphendemerits=0 \par\medskip\noindent}}
\newcommand{\si}{\sigma}
\newcommand{\G}{\mathcal G}
\newcommand{\cal}{\mathcal}
\newcommand{\sep}{^{\rm sep}}
\newcommand{\chara}{\mbox{\rm char}\,}
\newcommand{\R}{\mathbb R}
\newcommand{\Q}{\mathbb Q}
\newcommand{\N}{\mathbb N}
\newcommand{\Z}{\mathbb Z}
\newcommand{\F}{\mathbb F}
\newcommand{\C}{\mathbb C}
\newcommand{\cV}{\mathcal V}
\newcommand{\cF}{\mathcal F}
\newcommand{\cN}{\mathcal N}
\newcommand{\cM}{\mathcal M}
\newcommand{\cU}{\mathcal U}
\newcommand{\Fp}{\F_p}
\newcommand{\cL}{\mathcal{L}}
\newcommand{\res}{\textit{res}}
\newcommand{\cK}{\mathcal{K}}
\newcommand{\Div}{\text{div}}
\newcommand{\Def}{\mathrm{Def}}
\newcommand{\Dim}{\mathrm{dim}}
\newcommand{\gs}{\mathrm{gs}}
\newcommand{\dd}{\mathbf{d}}
\newcommand{\cc}{\mathbf{c}}
\newcommand{\llangle}{\langle\!\langle}
\newcommand{\rrangle}{\rangle\!\rangle}
\newcommand{\sS}{\mathbf{S}}
\newcommand{\ato}{\mathrm{at}}
\newcommand{\eq}{\mathrm{eq}}
\newcommand{\dcl}{\mathrm{dcl}}
\newcommand{\acl}{\mathrm{acl}}
\newcommand{\NIP}{\mathrm{NIP}}
\newcommand{\NTP}{\mathrm{NTP}}
\newcommand{\RCF}{\mathrm{RCF}}
\newcommand{\ACF}{\mathrm{ACF}}
\newcommand{\DCF}{\mathrm{DCF}}
\newcommand{\ACVF}{\mathrm{ACVF}}
\newcommand{\DOAG}{\mathrm{DOAG}}
\newcommand{\RCVF}{\mathrm{RCVF}}
\newcommand{\RV}{\mathrm{RV}}
\newcommand{\PCF}{p\mathrm{CF}}
\newcommand{\CODF}{\mathrm{CODF}}
\newcommand{\dom}{\mathrm{dom}}
\newcommand{\range}{\mathrm{range}}
\newcommand{\qf}{\mathrm{qf}}
\newcommand{\Th}{\mathrm{Th}}
\newcommand{\usd}{\mathrm{USD}}
\newcommand{\tame}{\mathrm{tame}}
\newcommand{\UD}{\mathrm{(UD)}}
\newcommand{\D}[1]{{\wideparen{#1}}}
\newcommand{\tp}{\mathrm{tp}}
\newcommand{\se}{\mathrm{se}}
\newcommand{\gse}{\mathrm{gse}}
\newcommand{\id}{\mathrm{id}}
\newcommand{\DD}{\mathbf{D}}
\newcommand{\SE}{\mathcal{S}\mathcal{E}}
\newcommand{\GSE}{\mathcal{G}\mathcal{S}\mathcal{E}}
\newcommand{\cO}{\mathcal{O}}
\newcommand{\J}{\nabla}
\newcommand{\cJ}{\mathcal{J}}
\newcommand{\cP}{\mathcal{P}}
\newcommand{\cA}{\mathcal{A}}
\newcommand{\cB}{\mathcal{B}}
\newcommand{\cC}{\mathcal{C}}
\newcommand{\cQ}{\mathcal{Q}}
\newcommand{\cI}{\mathcal{I}}
\newcommand{\cW}{\mathcal{W}}
\newcommand{\ord}{\text{ord}}
\newcommand{\Int}{\mathrm{Int}}
\newcommand{\cLr}{\mathcal{L}_\mathrm{ring}}
\newcommand{\In}{\text{Int}}
\newcommand{\ogr}{\mathrm{og}}
\newcommand{\CF}{\mathrm{CF}}
\newcommand{\cH}{\mathcal{H}}
\newcommand{\Pres}{\mathrm{Pres}}
\newcommand{\bard}[1]{\overline{{#1}}^\delta}
\newcommand{\bd}{\bar{\delta}}
\newcommand{\cZ}{\mathcal{Z}}
\newcommand{\cG}{\mathcal{G}}
\newcommand{\cS}{\mathcal{S}}
\newcommand{\da}{\underline{d}}

\begin{abstract}  We continue the study of a class of topological $\cL$-fields endowed with a generic derivation $\delta$, focussing on describing definable groups. We show that one can associate to an $\cL\cup\{\delta\}$ definable group a type $\cL$-definable topological group. We use the group configuration tool in o-minimal structures as developed by K. Peterzil. 
\end{abstract}

\thanks{$(\dagger)$ Research Director at the "Fonds de la Recherche
  Scientifique (F.R.S.-F.N.R.S.)"}

\maketitle
\section{Introduction}
\par Let $\cK$ be an $\cL$-structure expanding a field of characteristic $0$, endowed with a non-discrete definable field topology, as introduced by A. Pillay in \cite{P87}. We assume that $\cK$  is a model of an $\cL$-open theory $T$ of topological fields (see section \ref{sec:open} for a precise definition). The language $\cL$ is a (multisorted) language which is, on the field sort, a relational expansion of the ring language possibly with additional constants (with further assumptions on relation symbols, as defined in \cite{guzy-point2010}, \cite{CP}). In particular the theory $T$ is a complete $\cL$-theory admitting quantifier elimination on the field sort.

Examples of such theories $T$ are the theories of algebraically closed valued fields, real-closed fields, real-closed valued fields, p-adically closed fields or henselian valued fields of characteristic $0$. Note that the first four theories are not only dp-minimal but respectively $C$-minimal, o-minimal, weakly o-minimal, $p$-minimal.
\par W. Johnson established a link between topological fields and dp-minimal fields. He showed that if $\cK$ is an expansion of a field, $K$ infinite, and $dp$-minimal but not strongly minimal, then $\cK$ can be endowed with a non-discrete definable field topology, namely $\cK$ has a uniformly definable basis of neighbourhoods of zero compatible with the field operations \cite[Theorem 9.1.3]{JW}. Moreover, in models of the theory of $\cK$, the topological dimension coincides with the dp-rank.
This definable field topology is a V-topology and so it is induced either by a non-trivial valuation or an absolute value. 
\par Here given an $\cL$-open theory $T$ of topological fields, we consider the {\it generic} expansion of a model of $T$ with a derivation $\delta$. Namely, letting $\cL_\delta:=\cL\cup \{^{-1}\}\cup\{\delta\}$ and $T_{\delta}$ the $\cL_{\delta}$-theory consisting of $T$ together with the axioms expressing that $\delta$ is a derivation, we consider the class of existentially closed models of $T_\delta$. 
 When $T$ is a theory of henselian fields of characteristic $0$, we identify the class of existentially closed models of $T_{\delta}$ (see \cite[Corollary 3.3.4]{CP} for a precise statement): we give an explicit axiomatisation $T_\delta^*$ and we show various transfer of model-theoretic properties from the theory $T$ to the theory $T_\delta^*$. An easy but quite useful result is that $T_\delta^*$ admits quantifier elimination as well (on the field sort).
\par The main aim of the present paper is to show that given an $\cL_\delta$-definable group (on the field sort) in a model of $T_\delta^*$, one can associate a type $\cL$-definable group. 
There are two steps in the proof. First we give a more direct and slightly more informative and more general  proof that $T_\delta^*$ has the $\cL$-open core property \cite[Theorem 6.0.8]{CP}. Second we perform two constructions. Both were first done for o-minimal expansions of a group. The first construction is due to A. Pillay \cite{P88} who, on a definable group, put a definable topology for which the group operations become continuous and the other one is due to K. Peterzil  who, starting from a group configuration, constructed a type definable group. In order to perform those constructions, besides the $\cL$-open core property, one uses that the topological dimension function is a well-behaved dimension function in models of $T$ \cite[Proposition 2.4.1, Corollary 2.4.5]{CP} (the topological dimension coincides with the algebraic dimension, has the exchange property and is a fibered dimension function \cite{vdD89}).
Finally one uses the transfer of the elimination of imaginaries from $T$ to $T_\delta^*$ \cite[Theorem 4.0.5]{CP} under the $\cL$-open core property of $T_\delta^*$.  
In \cite{CP}, in the case the topology on the models of $T$ is given by a valuation, we need an intermediate result on continuity almost everywhere of definable functions to the value group and we only show it holds under certain conditions on the value group \cite[Proposition 2.6.11]{CP}. Here we proceed more directly constructing, what we call an $\cL$-definable envelop of an $\cL_\delta$-definable subset. It will help us to associate with an $\cL_\delta$-definable group a large definable $\cL$-definable subset of the so-called envelop on which we can recover the group operations on the differential points.

\

The contents of the paper are as follows. In section 2, we review the properties that we need on the topological dimension in models of an $\cL$-open theory; we recall the definition of the theory $T_\delta^*$ and various transfer properties between $T$ and $T_\delta^*$.
In section 3, we give a direct proof of the $\cL$-open core property of $T_\delta^*$ (Corollary \ref{opencore}): given a $\cL_\delta$-definable set $X$, one constructs a $\cL$-definable set that we call an $\cL$-open envelop  (Proposition \ref{prop:decomp2}). 
In section 4, we prove our main result, namely we show how to associate with an $\cL_{\delta}$-definable group a type $\cL$-definable topological group (Theorem \ref{local}). In section 5 (the annex), we revisit the question of when a topological differential field embeds in a model of the scheme (DL).

\

\noindent{\bf Acknowledgements} The first part of this paper owes a lot to a recent work with Pablo Cubid\`es on differential topological fiels; among other things, it helped the author to revisit a former work on definable groups. For the second part, the author is indebted to the work of Kobi Peterzil (and its nice exposition)  on group configuration in an o-minimal setting.

\section{Preliminaries}
\subsection{Model theory and topological fields} \label{sec:open}
Let  $\cL_{\mathrm{ring}}:=\{\cdot,+,-,0,1\}$ and $\cL_{\mathrm{field}}:=\cL_{\mathrm{ring}}\cup\{{\,}^{-1}\}$ denote the respective language of rings and of fields. Any field is an $\cL_{\mathrm{field}}$-structure by extending the multiplicative inverse to 0 by $0^{-1}=0$. 

We will follow standard model theoretic notation and terminology. Lower-case letters like $a, b, c$ and  $x,y,z$ will usually denote finite tuples and we let $\vert x\vert$ denote the length of $x$. We will sometimes use $\bar{x}$ to denote a tuple of the form $(x_1,\ldots, x_n)$ where each $x_i$ is itself a tuple. Given an $\cL$-structure $K$ and an $\cL$-formula $\varphi(x)$ with $\vert x\vert=n$, we let $\varphi(K)$ denote the set $\{a\in K^n\mid K\models \varphi(a)\}$.  Given a subset $X\subset K^n$ and $a\in K^{\ell}$, $0<\ell<n$, the fiber of $X$ over $a$ is denoted by $X_{a}:=\{b\in K^{n-\ell} \colon (a, b)\in X\}$. By an $\cL$-definable set, we mean definable with parameters. If we want to restrict the subset where the parameters vary we will use $\cL(A)$-definable, $A\subset K$. If we wish to specify that it is definable without parameters, we use $\cL(\emptyset)$-definable.

\

Let $\cL_r$ be a relational extension of $\cL_{\mathrm{field}}$ and assume that the language $\cL$ (possibly multi-sorted) extends $\cL_{r}$ in which every sort is an imaginary $\cL_r$-sort. In particular, the field-sort will always be the home sort of an $\cL$-structure. 
Let us recall the definition of $\cL$-open theories $T$. Let $K$ be a field of characteristic 0 in the language $\cL$ endowed with a definable field topology, namely there is  an $\cL$-formula $\chi(x,z)$ be  providing a basis of neighbourhoods of 0 \cite{P87}. Throughout the text, we will assume that the topology on $K$ is given by such formula $\chi$. Let $T$ be the $\cL$-theory of $K$. 

Then such theories $T$ with the following further assumptions are called $\cL$-open \cite{CP}:

\begin{enumerate}
\item[$(\mathbf{A})$] 
\begin{enumerate}
\item[(i)] $T$ has relative quantifier elimination with respect to the field sort, 
\item[(ii)] every quantifier-free $\cL$-definable subset of the field sort is a finite union of sets of the form: an intersection of an $\cL$-definable open set with a Zariski closed set.
\end{enumerate}
\end{enumerate}

\

One can choose the language $\cL$ in such a way that the following theories are examples of $\cL$-open theories: theories of real closed fields, $p$-adically closed valued fields, real closed valued fields, algebraically closed valued field of characteristic $0$ and 
the theories of any Hahn power series field $k(\!(t^\Gamma)\!)$, where $k$ is a field of characteristic 0 and $\Gamma$ is an ordered abelian group  \cite[Examples 2.2.1]{CP}, \cite[Proposition 4.3]{Fl}.
\subsection{Definable sets in models of $\cL$-open theories}

Let us recall a simple lemma in \cite{CP} on the form of $\cL$-definable subsets in models of $\cL$-open theories $T$. First we need to some notations.

\

Let $\cA$ be a finite subset of $K[x,y]$. We let 
$\cA^y:=\{P\in \cA\mid \deg_y(P)>0\}$, 
We let the $\cL_{\mathrm{ring}}(K)$-formula $Z_\cA(x,y)$ be $Z_\cA(x,y) := \bigwedge_{P\in \cA} P(x,y)=0$.  
Thus the algebraic subset of $K^{n+1}$ (Zariski closed set) defined by $\cA$ corresponds to $Z_\cA(K)$. For an element $R\in K[x,y]$ we let 
\begin{equation}\label{eq:locclosed}
Z_\cA^R(x,y) := Z_\cA(x,y)\wedge R(x,y)\neq 0. 
\end{equation}

\begin{lemma}\label{lem:niceform}\cite[Corollary ]{CP} Let $\cK$ be a model of an $\cL$-open theory $T$. Then every $\cL$-definable set $X\subseteq K^{n+1}$ (in the field sort) is defined by an $\cL$-formula $\varphi(x,y)$ with $x=(x_0,\ldots, x_{n-1})$ and $y$ a single variable such that
\[
\varphi(x,y) \leftrightarrow \bigvee_{j\in J} Z_{\cA_j}^{S_j}(x,y) \wedge \theta_j(x,y)
\]
with $J$ a finite set such that for each $j\in J$, $\theta_j$ is an $\cL$-formula that defines an open subset of $K^n$, $Z_{\cA_j}^{S_j}(x,y)$ is as defined in \ref{eq:locclosed} with either 
\begin{enumerate}
\item $\cA_j\subseteq K[x]$ and $S_j\in K[x,y]$ or 
\item  $\cA_j\subseteq K[x,y]$, $\cA^y=\{P\}$ and $S_j=(\frac{\partial}{\partial y} P) R_j$ for some $R_j\in K[x]$. 
\end{enumerate}
\end{lemma}

In a topological field $K$, one has a natural notion of topological dimension $\Dim$. For $X$ a definable subset of $K^n$, define $\Dim(X):=max\{\ell:\;$there is a projection $\pi^n_\ell:K^n\rightarrow K^{\ell}$ such that $\pi(X)$ has non-empty interior$\}$. Denote by $\overline{X}$ the closure of $X$ in the topological sense, let $\rm{Fr}(X):=\overline{X}\setminus X$. For $X\subset Y$, let $\rm{Int}_{Y}(X)$ be the subset of elements of $X$ which have a neighbourhood contained in $Y$ and $\rm{cl}_{Y}(X)$ be the subset of elements of $Y$ for which any neighbourhood has a non-empty intersection with $X$.

\ 

\subsection{Dimension functions} 

In this paragraph we fix a complete $\cL$-theory $T$ (without the previous assumptions on $\cL$), a sufficiently saturated model $\cU\models T$ endowed with a fibered dimension function $d$ (on definable sets) and we recall a few well-known properties of such fibered dimension. The notion of fibered dimension was introduced by L. van den Dries as follows. 
\begin{definition} \cite{vdD89} \label{fibdim}Let $\Def(\cU)$ be the set of $\cL$-definable subsets  $\cU$. A {\it fibered dimension function} $d: \Def(\cU)\rightarrow \{- \infty \}\cup On$ satisfies the following axioms. Let $S, S_{1}, S_{2}, T\in \Def(\cU)$.
\begin{enumerate}
\item $d(S)=-\infty$ iff $S=\emptyset$, $d(\{a\})=0$, for each $a\in U$.
\item $d(S_{1}\cup S_{2})=max\{d(S_{1}), d(S_{2})\}$.
\item Let $S\in U^m$, then for any permutation $\si$ of $\{1,\cdots,m\}$, we have $d(S^{\si})=d(S)$.
\item Let $S\subset U^{n+m}$, $S_{\bar x}:=\{\bar y\in A^m: (\bar x, \bar y)\in S\}$ and $S(\gamma):=\{\bar x\in U^n: d(S_{\bar x})=\gamma\},$ $\gamma\in On$. Then, $S(\gamma)\in Def(\cU)$ and $d(\{(\bar x,\bar y)\in S: \bar x\in S(\gamma)\})=d(S(\gamma))+\gamma.$
\end{enumerate}
If one adds that $d(U)=1$, one obtains a dimension function taking its values in $\{-\infty\}\cup \N$ and one can relax the condition $(4)$ by asking it only for $m=1$ (see \cite{vdD89} Proposition 1.4). 
\end{definition}
Given a (fibered) dimension function $d$ on $\cU$, one may extend it to the space of types $S_{n}^T(C)$ over a subset $C\subset U$ by
\[
d(p):=\inf\{d(\varphi(\cU,c): \varphi(x, y) \;{\rm an}\; \cL(\emptyset) {\rm -formula\; and}\; \varphi(x,c)\in p\}.
\]
Then we extend $d$ on a tuple $a$ of elements of $\cU$ by defining $d(a/C)$ as the dimension of the type $tp(a/C)$ of $a$ over $C$.
 One then can show \cite[Lemma 1.6]{BCP}:
Let $a, b$ be two tuples of elements of $U$, then: 
\begin{equation}\label{dim}
d(a b/C)=d(a/C\cup b)+d( b/C).
\end{equation}

Let $B\subset \cU$ and let $p(x)$ be a partial $n$-type. Let $X\subset U^n$ with $X=p(U)$. Then $X$ is {\it finitely satisfiable} in $B$ if for any formula $\varphi(x)\in p(x)$, we have $\varphi(M)\neq \emptyset.$

A tuple $a$ in a definable set $X$ is called {\it generic over $C$} (w.r. to the dimension d) if $d(X)=d(a/C)$.

\defn \label{dlarge} 
Let $B\subset A\subset U^n$ be two definable subsets of $\cU$, then $B$ is almost equal to $A$ (or large in $A$) if $d(A\setminus B)<d(A)$ \cite[section 2]{simon-walsberg2016}.
\edefn
Let us recall the following result which can be found in \cite{P88}, proven in the setting of o-minimal theories but it only uses the notion of a fibered dimension.
\fct \label{def-large} {\rm \cite[Proposition 1.13, Remark 1.14]{P88}} 

Let $A$ be an $\cL_{B}$-definable subset of $\cU$, with $B\subset U$ and let $\varphi( x; y)$ be an $\cL$-formula. Then $\{c:\;\varphi(U; c)\cap A$ is almost equal to $A\}$ is $B$-definable.

Moreover $\{c:\;\varphi(U; c)\cap A$ is almost equal to $A\}=\{c:$ for every $\cL$-generic point $u$ of $A$ over $c$, $\cU\models \varphi(u; c)\}$.
\efct

\subsection{A cell decomposition theorem}\label{sec:top}

Throughout this section, let $\cK\models T$, where $T$ is $\cL$-open. Let us first observe a few properties of the topological dimension $\Dim$ on definable subsets of $\cK$. First $\Dim$ is a fibered dimension function. The proof follows the same strategy as in \cite{vdD89} by noting that the dimension $\acl$-$\dim$ is induced by the field algebraic closure and has the exchange property (on the field sort) (namely $T$ is a geometric theory on the field sort) and that it coincides with the topological dimension (see \cite[Proposition 2.4.1]{CP}).

Note that we may also define the dimension of an $n$-tuple $\bar a$ of elements of $K$ as the cardinality of a maximal subtuple of $\acl$-independent elements and it coincides with the previous definition given above (see for instance \cite[section 4]{BCP}). 

Moreover, the dimension of the frontier of a definable set has the following property.
\begin{equation} \label{frontier} 
\Dim(\rm{Fr}(X))<\Dim(X)=\Dim(\overline{X}).
\end{equation}
As noted in \cite[Chapter 4, Corollary 1.9]{D98}, it implies for any definable set $X\subset Y$ that 
 \begin{equation}\label{relativefrontier}
 \dim(X\setminus \rm{Int}_{Y}(X))<\dim(Y).
 \end{equation}
Indeed, note that $X\setminus \rm{Int}_{Y}(X)=X\cap cl_{Y}(Y\setminus X)=cl_{Y}(Y\setminus X)\setminus(Y\setminus X)\subset Fr(Y\setminus X)$ and so $\dim(X\setminus \rm{Int}_{Y}(X))<\dim(Y)$.

\lmm \cite[Lemma 2.4]{P88} \label{coheir} 
Let $\tilde K$ be an elementary saturated extension of $\cK$ and $X$ be definable subset of $\tilde K$ with parameters in $K$. Let $a\in X$ be generic over $K$ and $c\in \tilde K$. Assume that $tp( a/K c)$ is finitely satisfiable in $\cK$, then $a$ is generic over $K c$.
\elmm
\pr It amounts to show that $\dim(a/K)=\dim(a/K c)$. Since $\dim(a \; c/K)=\dim(a/K c)+\dim(c/K)$ (see equation \ref{dim}), we will show that $\dim(a\;  c/K)=\dim(a/K)+\dim( c/K)$. 
Since $\dim=\acl$-$\dim$, let $a_1$, respectively $ c_1$, maximal $\acl$-independent subtuples of $a$, respectively $c$ over $K$. By the way of contradiction, suppose that $a_1 c_1$ is not $\acl$-independent over $K$, then w.l.o.g. we may assume that $c_{1}=c_{11}c_{12}$ with $\vert c_{11}\vert=1$ such that $c_{11}\in \acl(c_{12} a_1 K)$. 
Let $\varphi(x, z, y; u)$ with $u\subset K$ such that $\varphi(c_{11},c_{12}, a_1; u)$ holds and $\exists^{\leq n} x\;\varphi(x, c_{12}, a_1; u)$ for some natural number $n$.
Since $tp(a/K \;c)$ is finitely satisfiable in $\cK$, there is $b\subset K$ such that $\varphi(c_{11},c_{12}, b; u)\wedge \exists^{\leq n} x\;\varphi(x, c_{12}, b; u)$ holds. This contradicts that $c_1$ was chosen to be $\acl$-independent over $K$.
\qed

\

\par Finally, one has the following description of definable subsets of topological fields models of an $\cL$-open theory \cite{CP}; it is the analogue of the cell decomposition proven for dp-minimal fields (see \cite[Proposition 4.1]{simon-walsberg2016}). Before stating the result, we need to recall the notion of correspondences.
\begin{definition} \cite[section 3.1]{simon-walsberg2016}
A {\it correspondence} $f\colon E\rightrightarrows F$ consists of two definable subsets $E, F$ together with a definable subset $graph(f)$ of $E\times F$ such that
\[
0< \vert\{y\in F:\;(x,y)\in graph(f)\}\vert<\infty, {\rm for all}\;\; x\in E.
\]
The correspondence $f$ is {\it continuous} at $x\in E$ if for every $V\in \cB$ there is $U\in \cB$ such that $(f(x),f(x'))\in V$ whenever $(x,x')\in U$.
\end{definition}

A correspondence $f$ is an $m$-correspondence if for all $x\in E$,  $\vert\{y\in F:\;(x,y)\in graph(f)\}\vert=m$. We denote by $f(x)$ the set $\{y\in F:\;(x,y)\in graph(f)\}$. Note that a $1$-correspondence is a function (together with its domain and image). Moreover a continuous $m$-correspondence is locally given by $m$ continuous functions \cite[Lemma 3.1]{simon-walsberg2016} (see also \cite[Lemma 2.6.2]{CP}). A $m$-correspondence $f$ on an open definable set $U$ is continuous on an open subset of $U$ almost equal to $U$ \cite[Proposition 2.6.10]{CP}.

We follow the following convention: if $f\colon K^0\rightrightarrows K^{n}$, then $graph(f)$ is identified with a finite set and if $U$ is an open subset of $K^n$ and $f\colon U\rightrightarrows K^{0}$, $graph(f)$ is identified with $U$.

\begin{proposition} \cite[Theorem 2.7.1]{CP} \label{prop:cell} Let $T$ be an $\cL$-open theory of topological fields and $\cK$ be a model of $T$. Let $X$ be a definable subset of $K^n$. There are finitely many definable subsets $X_{i}$ with $X=\bigcup_i X_{i}$ such that $X_{i}$ is, up to permutation of coordinates,  the graph of a definable continuous $m$-correspondence $f\colon U_{i}\rightrightarrows K^{n-d}$, where $U_{i}$ is a definable open subset of $K^d$, for some $0\leq d\leq n$.
\end{proposition}
\subsection{Generic differential expansions of topological fields}\label{sec:dertop}
Let $T$ be an $\cL$-open theory of topological fields and $\cK$ be a model of $T$. Let $\cL_\delta$ be the language $\cL$ extended by a unary function symbol $\delta$. Denote by $\cK_\delta$ the expansion of $\cK$ to an $\cL_\delta$-structure. Let $T_\delta$ be the $\cL_{\delta}$-theory $T$ together with the usual axioms of a derivation, namely,  
\[
\begin{cases}
\forall x\forall y(\delta(x+y)=\delta(x)+\delta(y))\\
\forall x\forall y(\delta(xy)=\delta(x)y+x\delta(y)).
\end{cases}
\] 

\nota Let $K_\delta\models T_\delta$.
For~$m\geqslant 0$ and~$a\in K$, we define 
\[
\delta^m(a):=\underbrace{\delta\circ\cdots\circ\delta}_{m \text{ times}}(a), \text{ with $\delta^0(a):=a$,}
\]
and~$\bar{\delta}^m(a)$ as the finite sequence~$(\delta^0(a),\delta(a),\ldots,\delta^m(a))\in K^{m+1}$. 

Similarly, given an element~$a=(a_1,\ldots,a_n)\in K^n$, we will write~$\bar{\delta}^m(a)$ to denote the element~$(\bar{\delta}^m(a_1),\ldots,\bar{\delta}^m(a_n))\in K^{(m+1)n}$. 
Let $\bar m:=(m_{1},\ldots, m_{n})\in \N^n$ and $\vert \bar m\vert:=\sum_{i=1}^n m_{i}$. We will write~$\bar{\delta}^{\bar m}(a)$ to denote the element~$(\bar{\delta}^{m_{1}}(a_1),\ldots,\bar{\delta}^{m_{n}}(a_n))\in K^{\vert \bar m\vert+n}$. For notational clarity, we will sometimes use $\J_m$ instead of $\bd_m$, especially concerning the image of subsets of $K^n$. For example, when $A\subseteq K$, we will use the notation~$\J_m(A)$ for~$\{\bar{\delta}^m(a):\;a\in A\}$ instead of $\bd^m(A)$. Likewise for $A\subseteq K^n$, $\J_{\bar m}(A):=\{(\bar{\delta}^{m_{1}}(a_1),\ldots,\bar{\delta}^{m_{n}}(a_n)):\;a\in A\}\subseteq K^{\vert \bar m\vert+n}$. 

We will call $\bar{\delta}^{\bar m}(a)$ a differential tuple and we will sometimes use the notation $a^{\J}$, suppressing the index $\bar m$.
\enota

Given $x=(x_0,\ldots,x_n)$, we let $K\{x\}$ be the ring of differential polynomials in $n+1$ differential indeterminates $x_{0},\ldots, x_{n}$ over $K$, namely it is the ordinary polynomial ring in formal indeterminates $\delta^j(x_{i})$, $0\leq i\leq n$, $j\in \omega$, with the convention $\delta^0(x_{i}):=x_{i}$. We extend the derivation $\delta$ to $K\{x\}$ by setting $\delta(\delta^i(x_j))=\delta^{i+1}(x_j)$. By a rational differential function we simply mean a quotient of differential polynomials. 
 
For $P(x)\in K\{x\}$ and $0\leqslant i\leqslant n$, we let $\ord_{x_i}(P)$ denote the \emph{order of $P$ with respect to the variable $x_i$}, that is, the maximal integer $k$ such that $\delta^k(x_i)$ occurs in a non-trivial monomial of $P$ and $-1$ if no such $k$ exists. We let $\ord(P)$, the order of $P$, be $\max_{i} \ord_{x_i}(P)$. Suppose $\ord(P)=m$. For $\bar{x}=(\bar{x}_0,\ldots,\bar{x}_n)$ a tuple of variables with $\vert \bar{x}_i\vert=m+1$, we let $P^*\in K[\bar{x}]$ denote the corresponding ordinary polynomial such that $P(x)=P^*(\bar{\delta}^m(x))$. 

Suppose $\ord_{x_n}(P)=m\geqslant 0$. Then, there are (unique) differential polynomials $c_i\in K\{x\}$ such that $\ord_{x_n}(c_i)<m$ and 

\begin{equation}\label{eq:diffpol}
P(x)=\sum_{i=0}^d c_i(x)(\delta^m(x_n))^i. 
\end{equation}

The separant $s_{P}$ of $P$ is defined as 
\[s_{P}:=\frac{\partial}{\partial \delta^m(x_{n})}P\in K\{x\}.
\]
 We extend the notion of separant to arbitrary polynomials with an ordering on their variables in the natural way, namely, if $P\in K[x]$, the separant of $P$ corresponds to $s_P:=\frac{\partial}{\partial x_{n}}P\in K[x]$. By convention, we induce a total order on the variables $\delta^j(x_{i})$ by declaring that 
\[
\delta^k(x_{i})< \delta^{k'}(x_{j}) \Leftrightarrow 
\begin{cases}
i<j \\
i=j \text{ and } k<k'. 
\end{cases}  
\]
This order makes the notion of separant for differential polynomials compatible with the extended version for ordinary polynomials, \emph{i.e.}, 
$s_{P^*} = s_P^*$. 
 
 We define an operation on $K\{x\}$ sending $P\mapsto P^\delta$ as follows: for $P$ written as in \eqref{eq:diffpol}
\[
P(x)\mapsto P^\delta(x) = \sum_{i=0}^d \delta(c_i(x)) (\delta^m(x_n))^i. 
\] 
A simple calculation shows that 
\begin{equation}\label{eq:diffpol2}
\delta(P(x)) = P^\delta(x) + s_P(x)\delta^{m+1}(x_n).
\end{equation}

\ 

Now we will describe a scheme of $\cL_{\delta}$-axioms generalizing the axiomatization of closed ordered differential fields ($\CODF$) given by M. Singer in \cite{singer1978}.  Let $\chi(x,z)$ be an $\cL$-formula providing a basis of neighbourhoods of 0. For $a=(a_1,\ldots,a_n)$ with $\vert a_i\vert=\vert z\vert$, we let 
\[
W_{a}:=\chi(K,a_{1})\times \cdots\times \chi(K,a_{n}).
\]

\defn Set $T_{\delta}^*:=T_{\delta}\cup (\mathrm{DL})$, where $(\mathrm{DL})$ is the following list of axioms: for every differential polynomial $P(x)\in K\{x\}$ with $\vert x\vert=1$ and $\ord_x(P)=m$, for variables $u=(u_0,\ldots,u_m)$ with $\vert u_i\vert=\vert z\vert$ and $y=(y_0,\ldots,y_m)$
\[
\forall u \Big(\exists y(P^*(y)=0 \land s_P^*(y)\ne 0 )
\rightarrow \exists x\big(P(x)=0\land s_P(x)\ne 0\land
(\bd^m(x)-y)\in W_{u}\big)\Big).
\]
\edefn

As usual, by quantifying over coefficients, the axiom scheme $(\mathrm{DL})$ can be expressed in the language $\cL_\delta$. 

\

When $T=\RCF$, the theory $\RCF_\delta^*$ corresponds to $\CODF$ (which is consistent). When $T$ is either $\ACVF_{0,p}$, $\RCVF$, $\PCF_d$ or the $\RV$-theory of $\mathbb{C}(\!(t)\!)$ or $\mathbb{R}(\!(t)\!)$, the consistency of the theory $T_\delta^*$ follows by results in \cite[Corollary 3.8, Proposition 3.9]{guzy-point2010} (and see also \cite[Theorem 3.3.2]{CP}). 
 In \cite{guzy-point2010}, we showed how to embed a differential topological field $(K,\delta)$ which satisfied a hypothesis called (Hypothesis (I)) (see \cite[Definition 2.21]{guzy-point2010}), similar to largeness in this topological context into a differential extension model of the scheme (DL). 
Largeness is a notion introduced by Pop \cite{Pop}). In \cite{CP}, instead of assuming a largeness hypothesis, we worked with henselian topological fields. (Note that henselian fields are large \cite{Pop}.) In the annex, we will indicate the hypothesis we need to make these embedding proofs work (see section \ref{annex}).
 
\

An immediate consequence of the axiomatisation is the density property of differential points in open subsets of models of $T_{\delta}^*$.
\begin{lemma}[{\cite[Lemma 3.17]{guzy-point2010}}] \label{fact:density} Let $\cK_\delta\models T_{\delta}^{*}$. 
Let $O$ be an open subset of $K^n$. Then there is $a\in K$ such that $\bar \delta^{n-1}(a)\in O$. 
\end{lemma}
\pr Let $\bar u:=(u_{0},\cdots, u_{n-1})\in O$. We consider the differential polynomial $\delta^{n-1}(x)$. The corresponding algebraic polynomial is $x_{n-1}$.
Then we consider the differential equation $\delta^{n-1}(x)=u_{n-1}$; its separant is $1$ and so we can apply the scheme (DL) and get a differential solution close to the algebraic solution $\bar u$. So there is $a\in K$ such that $\delta^{n-1}(a)=u_{n-1}$ and $\bar \delta^{n-1}(a)\in O$.
\qed

\

Under assumption $(\mathbf{A})$ on the $\cL$-theory $T$, different model-theoretic properties transfer from $T$ to $T_\delta^*$, as shown by the following results. 

\begin{theorem}\label{thm:QE} \cite[Theorem 4.1]{guzy-point2010} \cite{CP} The theory $T_{\delta}^{*}$ admits quantifier elimination relative to the field sort in $\cL_{\delta}$. 
\end{theorem} 

Using this quantifier elimination result, one can easily show the following two transfer results. (A proof of the second result may be found in \cite[Appendix A.0.5]{CP}.)

\begin{proposition}[{\cite[Corollary 4.3]{guzy-point2010}}] The theory $T_\delta^*$ is $\NIP$, whenever $T$ is $\NIP$.
\end{proposition}

\begin{proposition}[Chernikov]\label{thm:cherni} The theory $T_\delta^*$ is distal, whenever $T$ is distal.
\end{proposition}

\

Let us recall the definition of open core (in a general setting) (see also \cite{DMS}).

\defn Let $\cK\models T$, let $\tilde \cL$ be an expansion of $\cL$ and let $\tilde T$ be the corresponding expansion of $T$. Let $\tilde \cK$ be an $\tilde \cL$-expansion of $\cK$. Then $\tilde \cK$  has $\cL$-open core if every $\tilde \cL$-definable open subset is $\cL$-definable.
An extension of $\tilde T$ has $\cL$-open core if every model of that extension has $\cL$-open core.
\edefn

Let $\cS$ be a collection of sorts of $\cL^{eq}$. We let $\cL^{\cS}$ denote the restriction of $\cL^{eq}$ to the home sort together with the new sorts in $\cS$.
\begin{theorem}\label{thm:EI} \cite[Theorem 4.0.5]{CP} Suppose that $T$ admits elimination of imaginaries in $\cL^\cS$ and that the theory $T_\delta^*$ has $\cL$-open core. Then the theory $T^*_{\delta}$ admits elimination of imaginaries in $\cL_{\delta}^{\cS}$. 
\end{theorem}
    
\par Finally let us recall a result on the existence of a fibered dimension function on definable subsets (on the home sort) in models of $T_\delta^*$ (which uses that $T_\delta^*$  admits quantifier elimination).
\defn An $\cL$-structure $\cM$ is called \textit{equationally bounded}
if  for each definable set $S\in M^{m+1}$ such that for every $\bar a\in M^m$, $S_{\bar a}$ is small, there exist finitely many  $\cL(M)$-terms $f_1(x_1,\ldots,x_m,y),\dotsc ,f_r (x_1,\ldots, x_m, y)$ such that for every $\bar a \in M^m$, there exists $1\leq i\leq r$ with $f_i(\bar a,y)\ne 0$ and $S_{\bar a}\subset \{b\in M:\,\,f_i(\bar a,b)=0\}$.
\edefn

\par Using a notion of independence built in from the algebra of all terms (t-independence), following \cite{vdD89}, \cite{Marc}, one can define on $\Def(\cK_\delta)$, when $\cK_\delta\models T_{\delta}^*$ a dimension function ($\Dim_{\delta}$)  \cite[Definition 2.3]{GP2}. In order to show that this dimension is fibered, one uses the closure operator $cl$ over $K$ which is defined by:  
$a\in cl(A)$ if and only if there is a differential polynomial $Q\in K \langle A\rangle\{X\}\setminus\{0\}$ such that $Q(a)=0$. 
Then associated with this closure operator, one has a notion of independence and dimension for tuples of elements. Let $\tilde \cK_\delta$ be a $\vert K\vert^+$-saturated extension of $\cK_\delta$.
Define for $\bar a=(a_1,\ldots,a_n)$ a tuple in $\tilde K$,
cl-$\dim(\bar a)$ := $\max \{ |B| : B\subseteq K\langle \bar a\rangle, B \;{\rm is\; cl-independent} \}.$  

Then one shows that, for $X\in \Def(\cK_\delta)$ (see \cite[Lemma 2.11]{GP2}):
\[
 \Dim_{\delta}(X)=\max\{{\rm cl-}\Dim(\bar c):\;c\in X(\tilde K)\}.
 \]
One then proves that any model of $T_{\delta}^*$ is equationally bounded, which entails that $\Dim_{\delta}$ is a fibered dimension function. Note that there are infinite definable subsets $X$ with $\Dim_{\delta}(X)=0$.
\begin{proposition}\label{cor delta dim}\cite[Corollary 3.10]{GP2}
Let $\cK_\delta\models T_{\delta}^*$, then there is a dimension function $\Dim_{\delta}$ that defines a fibered dimension function on $\Def(\cK_\delta)$.
\end{proposition}

This dimension function has been further investigated in \cite{BCP} and one can check that for $\bar a$ a tuple of elements in $\tilde K$
that $cl$-$\dim(\bar a$)=inf$\{\Dim_{\delta}(\varphi(K)):\; \varphi(\bar x)\in tp(\bar a/K)\}$. 

\section{The open core property in models of $T_{\delta}^*$}
From now on, let $\cK$ be a model of an $\cL$-open theory $T$ and let $\cK_\delta$ its expansion by a derivation $\delta$.  In this section to a $\cL_{\delta}$-definable set $X$, we will associate an $\cL$-definable set where the differential prolongation of $X$ is dense. Such result was already shown in \cite{CP} using a characterization of continuous $\cL_\delta$-correspondences (with $\cL$-definable domain).

\nota\label{nota:qf}
Under assumption $(\mathbf{A}(ii))$, any quantifier-free (relative to the field sort) $\cL_\delta$-definable set $X\subseteq K^n$ is of the form $\J_m^{-1}(Y)$ for a quantifier-free $\cL$-definable set $Y\subseteq K^{n(m+1)}$ (quantifier-free relative to the field sort). Indeed, let $x=(x_1,\ldots,x_n)$ be a tuple of field sort variables. By assumption on the language $\cL$, any $\cL_\delta$-term $t(x)$ is equivalent, modulo the theory of differential fields, to an $\cL_{\delta}$-term  $t^*(\delta^{m_1}(x_1),\ldots,\delta^{m_n}(x_n))$ where $t^*$ is an $\cL$-term, for some $(m_1,\ldots,m_n)\in \N^n$. Therefore, by possibly adding tautological conjunctions like~$\delta^k(x_i)=\delta^k(x_i)$ if needed, we may associate with any $\cL_\delta$-formula $\varphi(x)$ without field sort quantifiers, an equivalent $\cL_\delta$-formula (modulo the theory of differential fields) of the form $\varphi^{*}(\bar{\delta}^m(x))$ where $m\in \N$ and $\varphi^{*}$ is an $\cL$-formula without field sort quantifiers. The formula $\varphi^{*}$ arises by uniformly replacing every occurrence of~$\delta^m(x_i)$ by a new variable~$y_i^m$ in~$\varphi$ with the natural choice for the order of variables $\varphi^{*}(y_1^0,\ldots,y_1^m,\ldots,y_n^0,\ldots,y_n^m)$. Therefore, if $X$ is defined by $\varphi$, letting $Y$ be the set defined by $\varphi^*$ gives that $X=\J_m^{-1}(Y)$.  
\enota

\begin{definition}[Order]\label{def:order} Let $X\subseteq K^n$ be an $\cL_\delta$-definable set. Let $\bar d=(d_{1},\ldots,d_{n})\in \N^n$ and let $Z\subset K^{d_{1}+\ldots+d_{n}}$ such that  $X=\J_{\bar d}^{-1}(Z)$ 
The \emph{order of $X$}, denoted by $o(X)$, is the smallest integer $m$ such that  $m=\max_{1\leq i\leq n} d_{i}$ and $X=\J_{\bar d}^{-1}(Z)$. Let $\varphi$ be a field-sort quantifier-free $\cL_{\delta}$-formula such that $X=\{(a_{1},\ldots,a_{n})\in K^n: K\models \varphi^*(\bd^{d_{1}}(a_{1}),\ldots, \bd^{d_{n}}(a_{n}))\}$; we say that $d_{i}$ (also denoted by $d_{x_{i}}$) is the order of $x_{i}$ in $\varphi$, $1\leq i\leq n$.
\end{definition}

Note that $o(X)=0$ if and only if $X$ is $\cL$-definable. 

\defn\label{def:circ1} \cite[Definition 4.0.2]{CP} Let $X \subseteq K^n$ be a non-empty quantifier-free $\cL_\delta$-definable set (relative to the field sort). 
Let  $m$ be a positive integer and let $Z\subseteq K^{(m+1)n}$ be an $\cL$-definable set such that:
\begin{enumerate}
\item $x\in X$ if and only if $\J_{m}(x)\in Z$ and
\item $\overline{Z}=\overline{\J_{m}(X)}$.
\end{enumerate} 
Then we call $(X,Z,m)$ as above a linked triple.
\edefn
\begin{proposition}\label{thm:fermeture}\cite[Proposition 4.0.3]{CP} The theory $T_{\delta}^*$ has $\cL$-open core if and only if for every model $\cK_\delta$, for every $\cL_\delta$-definable set $X\subset K^n$, there is an integer $m$ and an $\cL$-definable set~$Z\subseteq K^{(m+1)n}$, such that $(X,Z,m)$ is a linked triple. In addition, if $T_{\delta}^*$ has $\cL$-open core, then there is a linked triple of the form $(X,Z,o(X))$.
\end{proposition}

Therefore, when $T_{\delta}^*$ has $\cL$-open core, we will associate two dimension functions to $X$, the first one $\Dim_{\delta}(X)$ (see Proposition \ref{cor delta dim}) and the second one $\Dim^*(X):=\min\{\frac{\Dim(Z)}{m+1}\colon (X,Z,m)$ is a linked triple$\}$.
\medskip
\par In \cite[Theorem 6.0.7]{CP}, we showed that $T_{\delta}^{*}$ has $\cL$-open core under following hypotheses on $T$. In the ordered case, the field sort is $\cL$-definably complete and in the valued case, the value group sort $\Gamma_{\infty}$ is $\cL$-definably complete. Furthermore, in the valued case, we assume the following on the value group: 
\begin{enumerate}
\item[$(\dagger)$] either there is a model $K'$ of $T$ for which $\Gamma(K')$ is a divisible ordered abelian group in which every infinite $\cL$-definable set has an accumulation point;
\item[$(\dagger\dagger)$] or there is a model $K'$ of $T$ with $\Gamma(K')=\Z$.  
\end{enumerate}
This last hypotheses enabled us to show that for every $\cL$-definable open set $V\subseteq K^n$, any $\cL$-definable function $f\colon V \to \Gamma_\infty$ is continuous almost everywhere \cite[Proposition 2.6.11]{CP}.

\subsection{Prolongations} 
At the beginning of this subsection, we simply assume that $\cK_\delta$ is a differential field of characteristic $0$. We introduce some notations on prolongations which will be used later on for  the $\cL$-open core property in models of $T_{\delta}^*$.

\begin{lemma}\label{lemdef:ratioprolong}\cite[Lemma 3.1.4]{CP} Let $x=(x_0,\ldots,x_n)$ be a tuple of variables and $y$ be a single variable. Let $P\in K\{x,y\}$ be a differential polynomial such that $k=\ord_y(P)\geqslant 0$. There is a sequence of rational differential functions $(f_i^P)_{i\geqslant 1}$ such that for every $a\in K^{n+1}$ and $b\in K$
\[
K\models [P(a,b)=0 \wedge s_P(a,b)\neq 0] \to \delta^{k+i}(b)=f_i^P(a,b).
\]
In addition, each $f_i^P$ is of the form 
\[
f_i^P(x,y) = \frac{Q_i(x,y)}{s_P(x,y)^{\ell_i}},
\]
where ${\ell_i}\in \N$, $\ord_{y}(Q_i)=\ord_y(P)$ and 
\[
\ord_{x_j}(Q_i)=
\begin{cases}
\ord_{x_j}(P)+i & \text{ if $\ord_{x_j}(P)\geqslant 0$ }\\
-1 & \text{ otherwise }
\end{cases}
\]
We call the sequence $(f_i^P)_{i\geqslant 1}$ the \emph{rational prolongation along $P$}. 
\end{lemma} 

\begin{notation}\label{not:prolong} Let $x=(x_0,\ldots,x_m)$ be a tuple of variables. For an integer $d\geqslant 0$, we define a new tuple of variables $x(d)$ which extends $x$ by $d$ new variables that is, 
\[
x(d):= (x_0, \ldots, x_m, x_{m+1},\ldots, x_{m+d}),
\]  
with the convention that if $d=0$, then $x(0)=x$,
When $\bar{x}=(\bar{x}_0,\ldots,\bar{x}_\ell)$ is a tuple of tuples of variables, we let $\bar{x}[d]:= (\bar{x}_0(d),\ldots,\bar{x}_\ell(d))$. Note that if $d\neq 0$ and $\bar{x}$ is not a singleton, then $\bar{x}[d]$ and $\bar{x}(d)$ are different.  

\

Let $x=(x_0,\ldots, x_n)$ and $y$ be a variable with $\vert y\vert=1$. Let $P\in K\{x,y\}$ be a differential polynomial of order $m$ and let $(f_i^P)_{i\geqslant 1}$ be its rational prolongation along $P$. 

 Let $\bar{y}$ be such that $\vert \bar{y}\vert=m+1$. We denote by $\lambda_{P}^d (\bar{x}[d], \bar{y},y_{m+1},\ldots,y_{m+d})$ the $\cL$-formula: 
\[
P^*(\bar{x}, \bar{y})=0 \wedge s_{P}^*(\bar{x}, \bar{y})\neq 0 \wedge \bigwedge_{i\geq 1}^d y_{m+i}=(f_i^P)^*(\bar{x}[d], \bar{y}).
\]
Given a tuple $a:=(a_{1},\ldots,a_{n})$ and an element $b$, we will call the tuple $(\bar{a}[d], \bar{b},b_{m+1},\ldots,b_{m+d})$ satisfying the formula $\lambda_{P}^d (\bar{a}[d], \bar{b},b_{m+1},\ldots,b_{m+d})$ the rational prolongation of $(a,b)$ (along $P$). When $d=0$, we set $\lambda_{P}^0 (\bar{x}[0], \bar{y})$ to be the formula $P^*(\bar{x}, \bar{y})=0 \wedge s_{P}^*(\bar{x}, \bar{y})\neq 0$.
\end{notation} 
We will make the following abuse of notation. When $P\in K[x,y]$ with $x=(x_0,\ldots,x_m)$ and $\vert y\vert=1$, we can view $P$ as $\tilde P^*$ where $\tilde P\in K\{x,y\}$ with $\ord_x(\tilde P)=m$. We will still denote by $\lambda_P^d(x(d),y(d))$ the formula $\lambda_{\tilde P}^d$.  

\

We will use the following notation concerning projections. 

\
\begin{notation}\label{sec:proj} For non-zero natural numbers $n,\;1\leqslant k\leqslant n$, we let $\pi_k\colon K^n \to K^k$ denote the projection onto the first $k$ coordinates and $\pi_{(k)}\colon K^n \to K$ denote the projection onto the $k^{th}$ coordinate. 
For natural numbers $n, m, \ell, k_i$ with $1\leq i\leq \ell\leq n$ and  $1\leqslant k_i\leqslant m+1$, we let $\pi_{[k_1,\ldots, k_{\ell}]}\colon K^{n(m+1)}\to K^{\ell(m+1)}$
denote the projection sending the $i^{\text{th}}$-block of $m+1$-coordinates to the first $k_i$-coordinates of such block, that is, 
\[
\pi_{[k_1,\ldots, k_{\ell}]}((x_{1,1},\ldots, x_{1,m+1},\ldots, x_{n,1},\ldots x_{n,m+1}))= (x_{1,1},\ldots, x_{1,k_1},\ldots, x_{\ell,1},\ldots x_{\ell,k_\ell}).
\] 
For $1\leq i\leq n$, we will allow in $\pi_{[k_1,\ldots, k_{n}]}$ the possibility for $k_{i}$ to be equal to $0$, in which case the corresponding subtuple $(x_{1,1},\ldots, x_{1,k_i})$ is simply empty.
Therefore the notation $\pi_{[0,\ldots, 0, k_{\ell}]}\colon K^{n(m+1)}\to K^{m+1}$ denotes the projection sending the $\ell^{\text{th}}$-block of $m+1$-coordinates to the first $k_\ell$-coordinates of such block, that is, 
\[
\pi_{[0,\ldots,0, k_{\ell}]}((x_{1,1},\ldots, x_{1,m+1},\ldots, x_{n,1},\ldots x_{n,m+1}))= (x_{\ell,1},\ldots x_{\ell,k_\ell}).
\] 
\end{notation}

Note that $\pi_{k_{1}}: K^{n(m+1)}\rightarrow K^{k_{1}}$ coincides with $\pi_{[k_{1}]}: K^{n(m+1)}\rightarrow K^{k_{1}}$, $k_{1}\leq m=1$.

\

The main result of this section is the following proposition, where we follow the same convention as in Proposition \ref{prop:cell}, namely an open subset of $K^0$ is a finite set of points.

\begin{proposition}\label{prop:decomp2} Let $\cK\models T_\delta^*$ and let $X\subseteq K^n$ be an $\cL_\delta$-definable set. Then there is $d\in \N$ and a finite family $\{Y_i\mid i\in I\}$ of $\cL$-definable subsets of $K^{n(d+1)}$ such that 
\[
X=\bigcup_{i\in I} \J_{d}^{-1}(Y_i), 
\]
and for each $i\in I$, $Y_i$ is either open, or equal to $\{\bd^d(a)\}$ for some $a\in X$, or the graph of a continuous $\cL$-definable correspondence $h_i$ such that  for some n-tuples $(d_1,\ldots,d_n)$ where $0\leq d_{i}\leq d+1$, $1\leq i\leq n$, $(\underline{d}_1,\ldots,\underline{ d}_n)$, where $\underline d_{i}:=\max\{d_{i},1\}$,
\begin{enumerate}
\item $\pi_{[\underline{d}_{1},\cdots,\underline{d}_{n}]}(Y_i)$ is open in $K^{d_{1}}\times\cdots\times K^{d_{n}}$,
\item  $h_i \colon \pi_{[\underline{d}_{1},\cdots,\underline{d}_{n}]}(Y_i) \rightrightarrows K^{d+1-d_{1}}\times\cdots\times K^{d+1-d_{n}}$, and
\item for every open subset $U\subseteq K^{nd}$ such that $U\cap \pi_{[\underline{d}_{1},\cdots,\underline{d}_{n}]}(Y_i)\neq \emptyset$, there is $a\in \pi_{[1,\ldots,1]}(U\cap \pi_{[\underline{d}_{1},\cdots,\underline{d}_{n}]}(Y_i))$ such that $\bd^d(a)\in Y_i$. 
\end{enumerate}  

We will denote $\bigcup_{i\in I} Y_i$ by $X^{**}$ and call it an $\cL$-open envelop of $X$ (in $K^{n(d+1)}$).
\end{proposition}
\begin{proof}
Let $m:=o(X)$; since $\J_{m}^{-1}$ commutes with finite union, by Theorem \ref{thm:QE}, Notation \ref{nota:qf} and Lemma \ref{prop:cell}, we may assume that $X$ is a finite union of subsets of the form $\J_m^{-1}(Z_{i})$, where $Z_{i}$ is an $\cL$-definable subset of  $K^{n(m+1)}$ which is either open, finite or, up to a permutation of coordinates, the graph of a continuous $\cL$-definable $\ell$-correspondence $f_i\colon U_{i} \rightrightarrows K^{n(m+1)-e}$, where $U_{i}\subset K^e$ is a non-empty open $\cL$-definable set, $1\leqslant e< n(m+1)$. We will consider each $Z_{i}$ separately. From now on, we drop the index $i$. The result is immediate for $Z$ in the following cases:
either $Z$ is open and we apply Lemma \ref{fact:density}, or $\dim(\pi_{[1,\ldots,1]}(Z))=0$ and so $\J_m^{-1}(Z)$ is finite. 

In the other cases, we proceed as follows. We associate with $Z$ a finite branching tree $(T,E)$ with root labelled by $Z=Z(0)$ and we denote the subset of the nodes at level $t$ by $T(t)$, $t\in \N$. If $v\in T(t)$, then we denote by $Z(v)$ the correspondence labelling the node $v$ with the property that $\nabla^{-1} Z(v)=\bigcup_{\substack{w\in T(t+1),\\ E(v,w)}} \nabla^{-1} Z(w)$ and $Z(v)\subset K^{n(m(v)+1)}$, where $m(v)$ has been chosen minimal such.
We will show that each branch is finite and at the end of each branch the correspondence which labels the vertex is one of the $Y_i$, $i\in I$ as in the statement of the proposition. Note that by Koenig's lemma, we get a finite tree. 
 
To simplify notation instead of denoting a correspondence obtained at level $t$ by $Z(v)$ with $v\in T(t)$, we will simply denote it by $Z(t)$ and we will denote $m(v)$ by $m(t)$.

Given a correspondence $Z(t)$, we associate a tuple $\bar d(t):=(d_{1}(t),\ldots,d_{n}(t))$ minimal in the lexicographic ordering, such that 
for each $1\leq i\leq n$, 
\[
\dim(\pi_{[0,\ldots,0,d_{i}(t)+1]}(Z(t)))=d_i(t),
\]
 placing ourselves in $K^{n(m(t)+2)}$  (see Notation \ref{sec:proj}).

For $d\in \N$, we use the notation $\da$ to mean $\max\{d,1\}$. On $\N^n$ we also put a partial ordering $\prec$ defined component-wise: $(d_1,\ldots,d_n)\prec (\tilde d_1,\ldots,\tilde d_n)$, if for all $1\leq i\leq n$, $d_i\leq \tilde d_i$ and for some $1\leq i\leq n$, $d_i<\tilde d_i$.

Define 
\begin{align*}
g(Z(t))&:=\bar d(t)\;\;{\rm  and}\\
 k(Z(t))&:=\vert\{i\colon 1\leq i\leq n, d_i(t)=m(t)+1\}.
 \end{align*}
  By assumption on $Z(0)$, we have $0\leq k(Z(0))<n$ and we will show that $k(Z(t))\leq k(Z(t-1))$, $t\geq 1$. 
We  decompose the tuple $\bar d(t)$ into two subtuples: $\bar d_{[k(Z(t))]}(t)$ of length $k(Z(t))$ which collects (respecting the order) all the components of $\bar d(t)$ equal to $m(t)+1$ and $\bar d_{[n-k(Z(t))]}$ such that $\bar d(t)=(\bar d_{[k(Z(t))]}, \bar d_{[n-k(Z(t))]})$. We proceed by induction on $(k(Z(t)),\bar d_{[n-k(Z(t))]}).$ Note that when $k(Z(t))=0$, then $\bar d_{[n-k(Z(t))]}=\bar d(t)=g(Z(t))$ and in that case we will show that when $t$ increases then $\bar d(t)$ decreases.
\par To the tuple $\bar d(t)$, we associate the following $n$-tuple of projections, with the convention that $d_i(t)^*+1:=\min\{d_i(t)+1,m(t)+1\}$, $1\leq i\leq n$, 
\[
\pi_{[d_1(t)^*+1]}, \pi_{[d_1(t)^*+1,d_2(t)^*+1]},\ldots, \pi_{[d_1(t)^*+1,\ldots,d_n(t)^*+1]}.
\]

Now let us describe how to get from $Z(t)$ to $Z(t+1)$. Since this step is similar to the step getting from $Z(0)\subset K^{n(m+1)}$ to $Z(1)$, in the discussion below, we will set $t=0$ and $d_i(0)=d_i$, $1\leq i\leq n$. 
 
If $d_{1}=m+1$, $\pi_{[m+1]}(Z)$ is an open subset of $K^{m+1}$, then we consider the least index $i$ such that $d_i\leq m$, $1<i\leq n$, and the set 
$\pi_{[m+1,\ldots,m+1, d_i+1]}(Z)$. (Note that by assumption on $Z$, there is such an index $i$.)

Or
$0\leq d_{1}<m+1$, so $\pi_{[d_{1}+1]}(Z)$ is a finite union of subsets $W_{1,j}$, $j\in J_1$ with $J_1$ finite, definable by $\cL$-formulas of the form 
\[ 
\varphi^1_{j}(x_1,y) : Z_{\cA_{1,j}}^{S_j}(x_1,y) \wedge \theta_{1,j}(x_1,y), \;\;(\dagger)_1
\]
with $\vert x_1\vert= d_{1}$, $\vert y\vert=1$, $\theta_j(x_1,y)$ an $\cL$-formula defining an open subset of $K^{d_{1}+1}$ and $\cA_{1,j}\subseteq K[x_1,y]$ a finite set of non-zero polynomials such that either $\cA_{1,j}\subseteq K[x_1]$ or $\cA_{1,j}^y=\{P_{1,j}\}$ and $S_j=(\frac{\partial}{\partial y} P_{1,j}) R_j$ for some $R_j\in K[x_1]\setminus \{0\}$ (see Lemma \ref{lem:niceform}). Note that every $\cA_{1,j}\neq\emptyset$ as otherwise we contradict the minimality of $d_{1}$. Moreover we may assume that for some formula $\varphi^1_{j}(x_1,y)$, we have $\cA_{1,j}\cap K[x_1]= \emptyset$  (by the assumption $\pi_{[d_{1}]}(Z)$ is a non-empty open subset unless $d_1=0$ and we use that $\dim$ coincides with the $\acl$-dimension (see section \ref{sec:top})). For each $j\in J$ such that $\cA_{1,j}\cap K[x_1]= \emptyset$, we may define the rational prolongation of $(x_1,y)$ along $P_{1,j}$ by the formula $\lambda_{P_{1,j}}^{\ell}(x_1,y(\ell))$, $\ell\geq 0$, where $y_{k}=f_{k}^{P_{1,j}}(x_1,y)$, $1\leq k\leq \ell$ and $y_0=y$ (see Notation \ref{not:prolong}).

In order to unify the notations, in case $d_1=m+1$, we still define $\lambda_{1}^{\ell}(x_1(\ell))$, $\ell\geq 0$ as the formula $x_1(\ell)=x_1(\ell)$ (namely putting no conditions on the prolongation of $x_1$).
We have $\lambda_{1}^{\ell}(K)=\pi_{[m+1]}(Z)\times K^{\ell}.$

In case $\cA_{1,j}\cap K[x_1]\neq \emptyset$, we consider the subset $\{\bar z\in Z\colon \pi_{[d_{1}]}(\bar z)\in W_j\}$ and we express it as a finite union of continuous correspondences $\tilde Z_j$, which will have the property that $\dim(\pi_{[d_1]}(\tilde Z_j))<d_1$.
So for each such correspondence $\tilde Z_j$, we have $g(\tilde Z_{j})\prec g(Z)$ and $d_{[n-k(\tilde Z_j)]}\prec d_{[n-k(Z)]}$. So we apply the induction hypothesis. So for each $\tilde Z_j$ we obtained a finite number of correspondences as described in the statement of the proposition and we label the corresponding number of nodes at level $1$ (that we connect to the root) by these correspondences. 

\

Either $d_{2}=m+1$ and so $\pi_{[0,d_{2}]}(Z)$ is an open subset of $K^{m+1}$. In this case we proceed, namely we consider

$\pi_{[d_1^*+1,m+1,d_{3}^*+1]}(Z)$, $0\leq d_{3}\leq m+1$.

Assume that $d_{2}<m+1$ and consider $\pi_{[d_1^*+1,d_{2}+1]}(Z)$. 
Again by Lemma \ref{lem:niceform}, we obtain that $\pi_{[d_1^*+1,d_{2}+1]}(Z)$ is a finite union of $\cL$-definable sets $W_{2,j}$, $j\in J_2$, each of which defined by a formula of the form:
\[ 
\varphi_j^{2}(x_{1},x_{2},y) := Z_{\cA_{2,j}}^{S_j}(x_{1},x_{2}, y) \wedge \theta_{2,j}(x_{1},x_{2}, y), \;\;(\dagger)_2
\]
with $\vert x_{1}\vert\leq d_1^*+1$, $\vert x_{2}\vert=d_{2}$, $\vert y\vert=1$, $\theta_{2,j}(x_{1},x_{2},y)$ an $\cL$-formula defining an open subset of $K^{d_1^*+1}\times K^{d_{2}+1}$ and $\cA_{2,j}\subseteq K[x_{1},x_{2},y]$ a finite non-empty set of non-zero polynomials such that either $\cA_{2,j}\subseteq K[x_{1},x_{2}]$ or $\cA_{2,j}^y=\{P_{2,j}\}$ and $S_2=(\frac{\partial}{\partial y} P_2) R_2$ for some $R_2\in K[x_{1},x_{2}]\setminus \{0\}$. (The fact that $\cA_{2,j}\neq\emptyset$ follows from the assumption on the dimension of $\pi_{[d_1^*+1,d_{2}+1]}(Z)$.) Moreover we may assume that for some formula $\varphi_j^{2}$ we have that $\cA_{2,j}\cap K[x_{1},x_{2}]= \emptyset$. For each such $j\in J_2$, we define the rational prolongation of $(x_1,x_2,y)$ along $P_{2,j}$. Namely let 
$\lambda_{P_{2,j}}^{\ell}(x_{1}(\ell),x_{2},y(\ell))$, $\ell\geq 0$, where $y_k=f_{k}^{P_{2,j}}(x_{1}(\ell),x_{2},y)$, $1\leq k\leq \ell$ and $y_0=y$ (see Notation \ref{not:prolong}). Note that letting $x_{1}=(x_{10},\ldots, x_{1m})$ and if $0\leq d\leq m$ is maximum $\deg_{x_{1d}} P_{2,j}>0$, then $x_{1}(\ell)=(x_{10},\ldots,x_{1(d+\ell)})$.

Again if $d_2=m+1$, we define $\lambda_{2}^{\ell}(x_{1},x_{2}(\ell))$, $\ell\geq 0$, as the formula $x_2(\ell)=x_2(\ell)$, putting no conditions on the prolongation of $(x_1,x_2,y)$.

 In case $\cA_{2,j}\cap K[x_{1},x_{2}]\neq \emptyset$, 
we consider the subset $\{\bar z\in Z\colon \pi_{[d_{1}^*+1,d_2]}(\bar z)\in W_{2,j}\}$ and we express it as a finite union of continuous correspondences $\tilde Z_{2,j}$, which will have the property that $\dim(\pi_{[d_1^*+1,d_2]}(\tilde Z_{2,j}))<d_1+d_2$.
So for each such correspondence $\tilde Z_{2,j}$, we have $g(\tilde Z_{2,j})\prec g(Z)$, $d_{[n-k(\tilde Z_{2,j})]}\prec d_{[n-k(Z)]}$ and we apply the induction hypothesis.

 Then we proceed with $\pi_{[d_1^*+1,d_2^*+1,d_{3}+1]}(Z)$ in case $d_3\leq m$, otherwise we proceed with $\pi_{[d_1^*+1,d_2^*+1,m+1, d_{4}+1]}(Z)$.
At step $1\leq k\leq n$ when $d_k\leq m$,  we denote the formulas we obtained by $\varphi^k_{j}$, $j\in J_k$ finite, and if $d_k=m+1$ we don't make any modifications 
but in order to uniformize notations we introduce  $\lambda_{k}^{\ell}(x_{1},\ldots,x_{k-1},x_{k}(\ell))$, $\ell\geq 0$, as the formula $x_k(\ell)=x_k(\ell)$, putting no conditions on the prolongation of $(x_{1},\ldots,x_{k-1},x_{k})$.
Then we consider the projection on the next coordinate.
 
Set for $1\leq i\leq n$:
\begin{enumerate}
\item in case $d_i\leq m$, $\ell_i:=m-d_i$,
\item in case  $d_i=m+1$, $\ell_i:=0$.
\end{enumerate} 
Then define $m(1)=m$ in case for all $i$, $d_i\leq m$, otherwise 
let 
\[m(1)=\max_{i}\{m+\sum_{j>i}^n \ell_j\colon 1\leq i<n, d_i=m+1\}.\]
 
 Then, in the construction we replace $Z$ by each of  the following subsets $\tilde Z$ as follows. Set $\bar z:=(z_1,\ldots,z_n)$ with $z_i:=(z_{i\,0},\ldots,z_{i\,m})$, and  for $1\leq i\leq n$,
 \par in case $d_i\leq m$, let $u_{i}:=(z_{i\,0},\ldots,z_{i\,d_i})$, 
\par in case $d_i=m+1$, let $u_i=z_i$.
\par For each $1\leq i\leq n$, we pick one of the formulas $\varphi^i_j(u_1,\ldots,u_i)$, $j\in J_i$, with $\vert u_i\vert=d_i^*+1$ and $\vert u_s\vert\leq d_s^*+1$, $1\leq s\leq i$, as occurring  above and we define 
 \[
\tilde Z:=\{\bar z[m(1)-m]:\ \bar z\in Z\wedge
 \bigwedge_{1\leq i\leq n} \lambda_{P_{i,j}}^{\ell_i+m(1)-m}(u_1(\ell_i+m(1)-m),\ldots, u_i(\ell_i+m(1)-m)\},
 \]

 making the convention that when $d_i=m+1$, the formula $\lambda_{P_{i,j}}^{\ell_i+m(1)-m}(u_1(\ell_i+m(1)-m),\ldots, u_i(\ell_i+m(1)-m))$ is replaced by $\lambda_{j}^{\ell_i+m(1)-m}(u_1,\ldots,u_{i-1}, u_i(\ell_i+m(1)-m))$ if $i>1$ and if $i=1$ by the formula $\lambda_{j}^{\ell_1+m(1)-m}(u_1(\ell_1+m(1)-m))$. Also note that if $d\leq d_i^*+1$ is maximal such that $\deg_{z_{s d}} P_{i,j}>0$, then in the formula $\lambda_{P_{i,j}}^{\ell_i+m(1)-m}(u_1(\ell_i+m(1)-m),\ldots, u_i(\ell_i+m(1)-m))$, $u_s(\ell_i)=(z_{s0},\ldots, z_{s d}, z_{s d+1}, \ldots, z_{s (d+\ell_{i}+m(1)-m}))$.
 
 Recall that along the way we took off $Z$ subsets of smaller dimension, so let us denote $Z_0$ the remaining subset. Then we get $\bigcup \nabla^{-1}(\tilde Z)=\nabla^{-1}(Z_0)$.

Then we apply Proposition \ref{prop:cell} to decompose each $\tilde Z$ into a finite union of correspondences to which we apply the preceding procedure. Let us denote one of these correspondences by $Z(1)$. 

Since there was nothing special going from $Z$ to $Z(1)$, replacing $0$ by $t$ (and $1$ by $t+1$), $m$ by $m(t)$, $d_i$ by $d_i(t)$, $\ell_i$ by $\ell_i(t)$, we have
from the construction above, that the tuple $\bar d(t)\in \N^n$ associated with $Z(t)$, $t\geq 1$, has the following properties:
\begin{enumerate}
\item if $d_i(0)<m+1$, then $d_i(1)\leq d_i(0)$, $1\leq i\leq n$, 
\item if $d_i(t)<d_i(t-1)$, then $d_i(t+1)\leq d_i(t)$, $1\leq i\leq n$,
\item if $d_i(t)>d_i(t-1)$, then $d_i(t-1)=m(t-1)$ and for some $i<j\leq n$, $d_{j}(t)<d_j(t-1)$, $1\leq i\leq n$,
\item $d_n(t)\leq d_n(t-1)\leq d_n(0)=d_n\leq m+1$,
\end{enumerate}
We make a similar abuse of notation by denoting $k(Z(t))$ by $k(t)$, and write the tuple $\bar d(t)$ as two subtuples one: $\bar d_{[k(t)]}(t)$ of length $k(t)$ which collects (in increasing order) all the components of $\bar d(t)$ equal to $m(t)+1$ and the other one $\bar d_{[n-k(t)]}$ such that $\bar d(t)=(\bar d_{[k(t)]}, \bar d_{[n-k(t)]})$. (Recall that $k(t+1)\leq k(t)<n$.)

So by $(4)$, for all $t$, $d_n(t)\leq d_n$, by $(1),\;(2)$ that the only indices $i$, $1\leq i\leq n$, for which $d_i(t+1)>d_i(t)$ are those where the projection onto the $i^{th}$-block of variables has the same dimension as the ambient space and by $(3)$ in order that this dimension increases from step $t$ to step $t+1$, the dimension of the projection onto some $j^{th}$-block of variables, $i<j\leq n$ has dimension strictly smaller than the dimension of the ambient space. 
Observe that if $k(t)=0$, then for all $i$, $d_i(t+1)\leq d_i(t)$ and 
$(d_1(t+1),\ldots,d_n(t+1))\prec (d_1(t),\ldots,d_n(t))\prec (d_1,\ldots,d_n)$, $t\geq 1$.

Furthermore if we go along a branch, then at some point $k(t)=k(t+1)$ and we cannot have $\bar d_{[n-k(t'+1)]}\prec \bar d_{[n-k(t')]}$ for infinitely many $t'>t$. So now suppose that $k(t)=k(t+1)$ and $\bar d_{[n-k(t+1)]}=\bar d_{[n-k(t)]}$. By $(3)$, we get that $\bar d_{[k(t+1)]}=\bar d_{[k(t)]}$. So $\bar d(t+1)=\bar d(t)$.

 Let us show that in this case, the differential points are dense in $Z(t+1)$ and that we can decompose $Z(t+1)$ into a finite union of correspondences ones for which we
 do have the required description and the other ones $\tilde Z(t+1)$ for which $g(\tilde Z(t+1))\prec g(Z(t))$. Let $\bar z:=(z_1,\ldots,z_n)\in Z(t+1)$ with $z_i:=(z_{i\,0},\ldots,z_{i\,m(t)})$, $1\leq i\leq n$.
When $d_i(t)\leq m(t)$, let $u_{i}:=(z_{i\,0},\ldots,z_{i\,d_i(t)})$ and $v_i:=(z_{i\,0},\ldots,z_{i\,d_i(t)-1})$. When $d_i(t)=m(t)+1$, $u_i=z_i=v_i$. 

We have $\pi_{[\bar d(t)]}(Z(t))=U_1\times\ldots\times U_n$ where $U_i\subset K^{d_i(t)}$ are open $\cL$-definable subsets of dimension $d_i(t)$ (with the convention that if $d_i(t)=0$, then $U_i$ is a point), $1\leq i\leq n$. By assumption $\pi_{[\bar d(t+1)]}(Z(t+1))=V_1\times\ldots\times V_n$ where $\emptyset\neq V_i\subset U_i$. 
  Let us describe a correspondence $f$ sending $(v_1,\ldots,v_n)\in V_1\times\ldots\times V_n$ to 
 $(u_1(\ell_1(t)),\ldots,u_n(\ell_n(t))$,
 where for each $1\leq i\leq n$, there is $j\in J_i$ such that $\lambda_{P_{i,j}}^{\ell}(u_1(\ell),\ldots,u_i(\ell))$ holds with the following convention:
 when $d_i(t)=m(t)+1$, then $\lambda_{P_{i,j}}^{\ell}(u_1(\ell),\ldots,u_i(\ell))$ is replaced by $\lambda_{i}^{\ell}(u_1(\ell),\ldots,u_{i-1}(\ell),u_i(\ell))$. Moreover we only keep those $(u_1(\ell_1(t)),\ldots,u_n(\ell_n(t)))$ which belongs to $Z(0)$. 

 The differential points are dense in $\pi_{[\bar d(t+1)]}(Z(t+1))$ (Lemma \ref{fact:density}), so we can pick a $n$-tuple of the form $(\delta^{d_1(t)-1}(a_1),\ldots, \delta^{d_n(t)-1}(a_n))$ close to $(v_1,\ldots,v_n)$. 
 Then we use the scheme $(DL)$ and the continuity of the rational functions $f_k^{P_{i,j}}$, $1\leq k\leq \ell$, associated with $P_{i,j}$, with $i$ such that  $d_i(t)\leq m(t)$ and $j\in J_i$.
We proceed by induction on $1\leq i\leq n$, namely we assume that the differential points are dense in $\pi_{[(i-1)(m(t)+1)]}(Z(t+1)$ with the convention that if $i=1$, this is the empty projection.
Now suppose that $d_i(t)\leq m(t)$, then for $(u_1,\ldots,u_{i-1})\in \pi_{[d_1(t)^*+1,\ldots,d_{i-1}(t)^*+1]}(Z(t+1))$, $u_i\in K^{d_i(t)+1}$, if 
\[
K\models P_{i,j}(u_1,\ldots,u_{i-1}, u_i)=0 \wedge s_{P_{i,j}}(u_1,\ldots,u_{i-1}, u_i)\neq 0,
\]
then there is a differential point $\bd^{d_i(t)+1}(u)$ close to $u_i$, for some $u\in K$, such that
\[
K\models P_{i,j}(u_1,\ldots,u_{i-1}, \bd^{d_i(t)+1}(u))=0 \wedge s_{P_{i,j}}(u_1,\ldots,u_{i-1}, \bd^{d_i(t)+1}(u))\neq 0.
\]

Let $U$ be a basic open neighbourhood of $(z_1,\ldots,z_i)\in \pi_{[i(m(t)+1)]}(Z(t+1)$. By induction hypothesis, there is a differential tuple close to $(z_1,\ldots,z_{i-1})$ and by scheme $(DL)$, there is a differential tuple close to $(z_1,\ldots,z_{i-1}, u_i)$. 
By the continuity of the functions $f_{k}^{P_{i,j}}$ (see Lemma \ref{lemdef:ratioprolong}), there is $V$ be a basic open neighbourhood of $(z_1,\ldots,z_{i-1},u_i)$ such that 

\[
V':=V\times f_1^{P_{i,j}}(V) \times \cdots \times f_{m(t)-d_i}^{P_{i,j}}(V) \subseteq U. 
\]

The rational functions $f_{k}^{P_{i,j}}$, $1\leq k\leq m(t)-d_i$,  applied to $\bd^{d_i(t)+1}(u)$ give its successive derivatives.

Then by \cite[Proposition 2.6.10]{CP}, the correspondence $f$ is continuous on an open subset $O$ almost equal to $V_1\times\ldots\times V_n$. It remains to apply induction to $f\restriction ((V_1\times\ldots\times V_n)\setminus O)$. 
\end{proof}

We apply the proposition above to obtain the open core property for $T_\delta^*$.

\kor \label{opencore} Let $K\models T_\delta^*$ and let $X\subseteq K^n$ be an $\cL_\delta$-definable set. Then there is an $\cL$-definable subset $Z$ such that $(X,Z,o(X))$ is a linked triple. In particular $T_{\delta}^*$ has $\cL$-open core.
\ekor
\pr By the above proposition, given an $\cL_{\delta}$-definable set $X$, there is an associated linked triple $(X,X^{**},d)$ with $X^{**}$ an open $\cL$-envelop for $X$ and $d\in \N$. In general this integer $d$ will be larger than $o(X)$ but, by Proposition \ref{thm:fermeture}, one can construct another linked triple with $d=o(X)$.
\qed

\medskip
\par The following corollary slightly improves \cite[Proposition 6.1.1]{CP} where we showed that one can associate with an $\cL_\delta$-definable $\ell$--correspondence, an $\cL$-definable $\ell$-correspondence. There we assumed the domain of the correspondence to be $\cL$-definable (but the same proof works when it is only $\cL_\delta$-definable). Furthermore, here we want the domain of the correspondence to have the special form described in Proposition \ref{prop:decomp2}. 
One can check that indeed the previous proof easily adapts.

\kor\label{prop:honest2} Let $K\models T_\delta^*$. Let $X\subset K^n$ and assume it is $\cL_{\delta}$-definable.
Let  $f\colon X\rightrightarrows K^{d}$ be an $\cL_\delta$-definable $\ell$-correspondence with $d, \ell\geq 1$. Let $X^{**}$ be an $\cL$-open envelop of $X$. There is $m\in \N$, 
an $\cL$-definable $\ell$-correspondence $F\colon X^{**}\rightrightarrows K^{d}$ such that for every $x\in X$ 
\[
f(x)=F(\bar{\delta}^m(x)).
\]
\qed
\ekor

\section{Definable groups in models of $T_\delta$}

First we will recall a few facts about definable groups and generics in the setting of o-minimal theories \cite{P88}, \cite{Pet},
since most of these notions only uses properties of the dimension function that hold in our present setting.
\par Then given an $\cL_\delta$-definable group $G$ in a model of $T_\delta^*$, where $T$ is an $\cL$-open theory of topological fields, we will show that on large $\cL$-definable subset of an $\cL$-open envelop of $G$ as defined in Proposition \ref{prop:decomp2}, one can recover an $\cL$-definable operation induced by the $\cL_\delta$-definable group law of $G$. 
\par Finally we will use the work of K. Peterzil on the group configuration in o-minimal structures \cite{Pet}. 
In a complete o-minimal theory admitting elimination of imaginaries, he showed that a group configuration gives rise to a transitive action of a type-definable group on infinitesimal neighbourghoods. Without appealing directly to a group configuration but using the same strategy, we will associate to an $\cL_\delta$-definable group a type $\cL$-definable group (possibly in a higher dimensional space).

In the last part of this section we will add sorts $\cS$ of  $\cL^{\mathrm{eq}}$ to the language $\cL$ in order that $T$ admits elimination of imaginaries in $\cL^{\cS}$, the restriction of $\cL^{\eq}$ to $\cL$ together with $\cS$. Note that on the one hand, the expansion $T^{\cS}$ of the theory $T$  in $\cL^{\cS}$ is again an $\cL^{\cS}$-open theory of topological fields (see \cite[Remark 2.2.1]{CP}) and on the other hand, by Corollary \ref{opencore} and \cite[Theorem 4.0.5]{CP}, $T_{\delta}^{\cS}$ admits elimination of imaginaries. 

\

We will first place ourselves in any $\cL$-structure $\cM$ (without our previous assumptions on the language $\cL$) but we will assume that $\cM$ is endowed with a fibered dimension which has the exchange property and that this dimension is preserved under definable bijection.

\par By definable group $\cG:=(G,\cdot,1)$ in $\cM$, we mean that the domain of $G$ is an $\cL$-definable subset of some cartesian product $M^n$ of $M$ and the graph of the group operation $\cdot$ is an $\cL$-definable subset of $M^{3 n}$. 

\fct \label{generic} {\rm \cite{P88}} Let $\cG$ be a definable group in $\cM$. Any element of $G$ is the product of two generic elements and given a generic element $a$ of $G$ over $b$, then $b\cdot a$ is generic in $G$.
\efct
\defn {\rm \cite{P88}} Let $\cG$ be a definable group in $\cM$. A {\it generic} definable subset of $G$ is a subset of $G$ such that finitely many translates cover $G$. 
\edefn
Note that it entails that a generic set has the same dimension than $G$ and so it contains a generic point of $G$.
\fct {\rm \cite[Lemma 2.4]{P88}} Let $\cG$ be a definable group of $\cM$. Let $X$ be a definable subset of $G$, almost equal to $G$, then $X$ is  generic. 
\efct

\

\par For the rest of this section, we will work in models of an $\cL$-open theory $T$ of topological fields. In particular, we assume that $T$ satisfies hypothesis ${\bf (A)}$.

\prop \label{largeV} Let $\cK\models T_{\delta}^*$ be sufficiently saturated, let $\cG:=(G,f_{\circ} ,f_{-1},e )$ be an $\cL_{\delta}$-definable group in $\cK$ (possibly with parameters) and let
$G^{**}$ be an $\cL$-open envelop of $G$.
Then there exist $\cL$-definable maps $F_{-1}:G^{**}\to G^{**}$ and $F_{\circ}:G^{**}\times G^{**}\to G^{**}$, a large open subset $V$ of $G^{**}$ and $Y$ a definable large open subset of $G^{**}\times G^{**}$ 
such that 
\begin{enumerate}
\item the map $F_{-1}:G^{**}\to G^{**}$ (respectively $F_{\circ}:G^{**}\times G^{**}\to G^{**}$) coincides on differential tuples with $f_{-1}$ (respectively $f_{\circ}$),
\item the map $F_{-1}: V\to V$ is a continuous idempotent map,
\item the map $F_{\circ}: Y\to V$ is continuous, 
\item for any $a\in V$, if $b$ is generic of $V$ over $a$, then $(b,a)\in Y$ and $(F_{-1}(b),F_{\circ}(b,a))\in Y$.
\end{enumerate}
\eprop
\pr The proof follows the same pattern as in \cite[Proposition 2.5]{P88} and we will construct a large subset $V$ of $G^{**}$ with the required properties by steps.

Assume that the domain of $\cG$ is included in some $K^{n}$ and is an $\cL_\delta$-definable set. Then by Proposition \ref{prop:decomp2}, there is an $\cL$-open envelop $G^{**}$ of $G$, described as follows. For some $d\in \N$, there is a finite family $\{Y_i\mid i\in I\}$ of $\cL$-definable subsets of $K^{n.(d+1)}$ such that $G^{**}=\bigcup_{i\in I} Y_i$, 
\[
G=\bigcup_{i\in I} \J_{d}^{-1}(Y_i), 
\]
and for each $i\in I$, $Y_i$ is either open, or equal to $\{\bd^d(a)\}$ for some $a\in G$, or the graph of a continuous $\cL$-definable correspondence $h_i$ such that  for some n-tuple $(d_1,\ldots,d_n)$:

\begin{enumerate}
\item $1\leqslant d_j<d+1$, $1\leq j\leq n$,
\item $\pi_{[d_{1},\cdots,d_{n}]}(Y_i)$ is open in $K^{d_{1}}\times\cdots\times K^{d_{n}}$,
\item  $h_i \colon \pi_{[d_{1},\cdots,d_{n}]}(Y_i) \rightrightarrows K^{d+1-d_{1}}\times\cdots\times K^{d+1-d_{n}}$, and
\item for every open subset $U\subseteq K^{nd}$ such that $U\cap \pi_{[d_{1},\cdots,d_{n}]}(Y_i)\neq \emptyset$, there is $a\in \pi_{[1,\ldots,1]}(U\cap \pi_{[d_{1},\cdots,d_{n}]}(Y_i))$ such that $\bd^d(a)\in Y_i$. 
\end{enumerate}  
Moreover letting $f_{-1}$ the inverse function on $G$, by Corollary \ref{prop:honest2}, there is a definable $\cL$-function $F_{-1}$ on $G^{**}$ which coincides with $f_{-1}$ on differential points. Similarly given $f_{\circ}$ the group law on $G\times G$ there is a definable $\cL$-function $F_{\circ}$ on $G^{**}\times G^{**}$ which coincides with $f_{\circ}$ on differential points. (We also use here that a $1$-correspondence is a function).
\par Let $I_{0}:=\{i\in I\colon \dim(Y_{i})=\dim(G^{**})\}$. By \cite[Proposition 2.6.10]{CP}, $F_{-1}\circ h_{i}$ is continuous on a large open subset $U_{i}$ of $\pi_{[d_{1},\cdots,d_{n}]}(Y_i)$. 

Let $\tilde Y_{i}:=h_{i}(U_{i})$ for $i\in I_{0}$ and 
\[
V_{1}:=\bigcup_{i\in I_{0}} \tilde Y_{i}.
\]
\cl \label{$V_{1}$} Let $V_{1}':=\{\bar x\in V_{1}: F_{-1}(F_{-1}(\bar x))=\bar x\}$. Then $V_{1}'$ is a large definable subset of $V_{1}$.
\ecl
\prcl Suppose on the contrary that we can find an open subset $U\subset V_{1}$ where  $F_{-1}(F_{-1}(\bar x))\neq \bar x$. Then choose in $U$ an element $\bar u^{\J}$ for some $\bar u\in G$. Since $F_{-1}$ and $f_{-1}$ coincide on differential points, we get a contradiction. 
\qed

\par From now on, we allow ourselves to replace $V_{1}$ by $V_{1}'$, namely we will assume that $F_{-1}\circ F_{-1}$ is the identity on $V_{1}$.

\

\par Using again \cite[Proposition 2.6.10]{CP}, $F_{\circ}\circ (h_{i},h_{j})$ is continuous on a large open subset $O_{i,j}$ of $\pi_{[d_{1},\cdots,d_{n}]}(Y_i)\times \pi_{[d_{1},\cdots,d_{n}]}(Y_j)$. 
Let $\tilde Y_{i,j}:=(h_{i},h_{j})(O_{i,j})$ for $i, j\in I_{0}$ and 
\[
Y_{0}:=\bigcup_{i, j\in I_{0}} \tilde Y_{i,j}.
\]
\cl \label{$V_{0}$} Let $V_{0}:=\{\bar y\in V_{1}: \forall \bar a\in G^{**}$ $\cL$-generic of $G^{**}$ over $\bar y$, $(\bar a,\bar y)\in Y_{0}\}$. Then $V_{0}$ is definable and almost equal to $V_{1}$. 
\ecl
\prcl The set $V_{0}$ is $\cL$-definable by Fact \ref{def-large}.
By the way of contradiction, assume there is a relatively open subset $U_{1}$ of $V_{1}\subset G^{**}$ such that if $\bar y\in U_{1}$, we have that $\{\bar a\in G^{**}:\;(\bar a,\bar y)\in Y_{0}\}$ is not almost equal to $G^{**}$. 
In particular, we may choose in $U_{1}$ a generic element $\bar y_{0}\in V_{1}$ of $G^{**}$ such that the set $\{\bar a\in G^{**}:\;(\bar a,\bar y_{0})\in Y_{0}\}$ is not almost equal to $G^{**}$. 
So its complement would contain an open subset $U_{2}(\bar y_{0})$.  Choose $\bar b_{0}\in G^{*}$ generic over $\bar y_0$ in that open set. Then $(\bar b_0,\bar y_0)$ is generic in $G^{**}\times G^{**}$ by equation \ref{dim} but would not belong to $Y_{0}$, a contradiction. 
\qed

\cl \label{ageneric} Let $a, b\in G$ and choose $a^{\J}$ $\cL$-generic over $b^{\J}$. Then $F_{\circ}(a^{\J},b^{\J})$ is $\cL$-generic over $b^{\J}$ and so it belongs to $V_{0}$.
\ecl
\prcl By hypothesis, $\dim(G^{**})=\dim(a^{\J}/b^{\J})$. Let us show that $\dim(F_{\circ}(a^{\J},b^{\J})/b^{\J})=\dim(a^{\J}/b^{\J})$.
By construction, we have that $F_{\circ}$ and $F_{-1}$ coincide respectively on $G^{\J}$ with $f_{\circ}$ and $f_{-1}$ respectively. So $F_{\circ}(a^{\J},b^{\J})=f_{\circ}(a,b)^{\J}$, $F_{-1}(b^{\J})=f_{-1}(b)^{\J}$ and $F_{\circ}(F_{\circ}(a^{\J},b^{\J}),F_{-1}(b^{\J}))=a^{\J}$.

\noindent So, $acl(b^{\J}, F_{\circ}(a^{\J},b^{\J}))\subset acl(b^{\J}, a^{\J})\subset acl(F_{\circ}(a^{\J},b^{\J}),F_{-1}(b^{\J}),b^{\J})\subset acl(F_{\circ}(a^{\J},b^{\J}),b^{\J})$. \qed

\cl \label{aigeneric} Let $a, b\in G$ and choose $a^{\J}$ $\cL$-generic over $b^{\J}$. Then $F_{-1}(a^{\J})$ is $\cL$-generic over $b^{\J}$ and so it belongs to $V_{1}$.
\ecl
\prcl By hypothesis, $\dim(G^{**})=\dim(a^{\J}/b^{\J})$. Let us show that $\dim(F_{-1}(a^{\J})/b^{\J})=\dim(a^{\J}/b^{\J})$.
By construction, we have that $F_{-1}$ coincide on $G^{\J}$ with $f_{-1}$. So $F_{-1}(a^{\J})=f_{-1}(a)^{\J}$ and $F_{-1}(F_{-1}(a^{\J}))=F_{-1}(f_{-1}(a)^{\J})=f_{-1}(f_{-1}(a)^{\J})=a^{\J}$. So, $acl(b^{\J}, a^{\J})\subset acl(b^{\J}, f_{-1}(a)^{\J})\subset acl(b^{\J}),a^{\J})$. \qed

\begin{align*}
{\rm Let}\;\;\tilde Y_{0}(\bar z):=\{(\bar x,\bar y)\in Y_{0}:\;&(\bar y,\bar z)\in Y_{0}\;\&\\
&F_\circ(\bar y,\bar z)\in V_{0}\;\&\\
&(\bar x,F_\circ(\bar y,\bar z))\in Y_{0}\;\&\\
&F_\circ(\bar x,F_\circ(\bar y,\bar z))=F_\circ(F_\circ(\bar x,\bar y),\bar z)\}
\end{align*}
\begin{align*}
{\rm Let\;} V_{0}^{''}:=\{\bar z\in V_{0}:\;\tilde Y_{0}(\bar z)\; {\rm is\; almost\;equal\;to\;} Y_{0}\}.
\end{align*}
 
The set $V_{0}^{''}$ is definable since being almost equal is an $\cL$-definable property and all the other data is $\cL$-definable.

\cl \label{$V_{0}^{''}$}$V_{0}^{''}$ is almost equal to $V_{0}$. 
\ecl
\prcl We proceed by contradiction. If not there would exist $\bar z\in V_{0}$ and an open neighbourhood $U$ of $\bar z$ such that for any element $\bar u\in U$, the set $\tilde Y_{0}(\bar u)$
is not almost equal to $Y_{0}$, which means that there exists an open subset $W$ of $Y_{0}$ containing $(\bar x,\bar y)$ where one of the following statement fails: 
 $(\bar y,\bar u)\in Y_{0},\;F_\circ(\bar y,\bar u)\in V_{0},\;(\bar x,F_\circ(\bar y,\bar u))\in Y_{0}$, 
 or if everything else hold, that $F_\circ(\bar x,F_\circ(\bar y,\bar u))=F_\circ(F_\circ(\bar x,\bar y),\bar u)$ fails.
 \par We may choose such $\bar u\in V_{0}\cap U$ of the form $\bar u^{\J}$ for some $\bar u\in G$ ($G^{\J}$ is dense in $G^{**}$). We choose $(\bar a^{\J},\bar b^{\J})\in W$ as follows.
 First we choose $\bar b^{\J}$ $\cL$-generic over $\bar u^{\J}$ (in particular $\bar b^{\J}\in V_1$). By Claim \ref{$V_{0}$}, $(\bar b^{\J},\bar u^{\J})\in Y_{0}$ and 
  by Claim \ref{ageneric}, $F_\circ(\bar b^{\J},\bar u^{\J})\in V_{0}$. (Note that $f_{\circ}(b,u)^{\J}=F_\circ(\bar b^{\J},\bar u^{\J}).$)
 Then we choose $\bar a^{\J}$ $\cL$-generic over $\bar b^{\J}$ and over $f_{\circ}(b,u)^{\J}$.  So $(\bar a^{\J},\bar b^{\J})\in Y_0$ and $(\bar a^{\J}, f_{\circ}(b,u)^{\J})\in Y_{0}$. Since $G$ is a group and $F_{\circ}$ and $f_{\circ}$ coincide on differential points, we get that $F_{\circ}(\bar a^{\J},F_\circ(\bar b^{\J},\bar u^{\J}))=F_\circ(F_\circ(\bar a^{\J},\bar b^{\J}),\bar u^{\J})$. 
 \qed

\begin{align*}
{\rm Let}\;\;V_{1}^{''}:=\{\bar z\in V_{0}:\{ \bar x\in V_{1}:\;&(F_{-1}(\bar x),\bar z)\in Y_{0}\;\&\\
&\;(\bar x,F_{\circ}(F_{-1}(\bar x),\bar z))\in Y_{0}\;\&\\
&\;F_{\circ}(\bar x,F_{\circ}(F_{-1}(\bar x),\bar z))=\bar z\}\; {\rm is \;\;almost\; equal \;to \;}V_{0}\}.
\end{align*}
The set $V_{1}^{''}$ is definable since being almost equal is an $\cL$-definable property.
\cl \label{$V_{1}^{''}$} $V_{1}^{''}$ is almost equal to $V_{0}$. 
\ecl
\prcl Suppose not then there would exist an open subset $U$ of $V_{0}$ over which  $\{\bar x\in V_{1}:(F_{-1}(\bar x),\bar z)\in Y_{0}\;\&
\;(\bar x,F_{\circ}(F_{-1}(\bar x),\bar z))\in Y_{0}\;\&\;F_{\circ}(\bar x,F_{\circ}(F_{-1}(\bar x),\bar z))=\bar z\}$ is not almost equal to $V_{0}$.

So we can find an element of the form $\bar z^{\J}\in U$, $\bar z\in G$ and $\bar x^{\J}\in V_{0}$ generic over $\bar z^{\J}$. By Claim \ref{aigeneric},
$F_{-1}(x^{\J})$ is generic  over $\bar z^{\J}$ and so it belongs to $V_{0}$ and $(F_{-1}(\bar x^{\J}),\bar z^{\J})\in Y_{0}$. 
By Claim \ref{ageneric}, $F_{\circ}(F_{-1}(\bar x^{\J}),\bar z^{\J})$ is generic over $\bar z^{\J}$.
Now let us show that $\bar x^{\J}$ is generic over $F_{\circ}(F_{-1}(\bar x^{\J}),\bar z^{\J})$. 
\par \noindent We first note that $acl(\bar z^{\J},\bar x^{\J})=acl(\bar x^{\J},F_{\circ}(F_{-1}(\bar x^{\J}),\bar z^{\J}))$. 
\par \noindent So $\dim(\bar x^{\J}, \bar z^{\J})=\dim(\bar x^{\J},F_{\circ}(F_{-1}(\bar x^{\J}),\bar z^{\J})).$ By equation \ref{dim}, 
\begin{align*}
\dim(\bar x^{\J}, \bar z^{\J})=&\dim(\bar x^{\J}/\bar z^{\J})+\dim(\bar z^{\J}) {\;\rm and \;}\\
\dim(\bar x^{\J},F_{\circ}(F_{-1}(\bar x^{\J}),\bar z^{\J}))=&\dim(\bar x^{\J}/F_{\circ}(F_{-1}(\bar x^{\J}),\bar z^{\J}))+\dim(F_{\circ}(F_{-1}(\bar x^{\J}),\bar z^{\J})).
\end{align*}
So, $\dim(\bar x^{\J}/F_{\circ}(F_{-1}(\bar x^{\J}),\bar z^{\J}))=\dim(\bar x^{J}/\bar z^{\J})$, so $\bar x^{\J}$ is generic over $F_{\circ}(F_{-1}(\bar x^{\J}),\bar z^{\J})$ and therefore 
$(\bar x^{\J}, F_{\circ}(F_{-1}(\bar x^{\J}),\bar z))\in Y_{0}$. Since $F_{\circ}$ and $F_{-1}$ coincide with $f_{\circ}$ and $f_{-1}$ on differential points, we get $F_{\circ}(\bar x^{\J},F_{\circ}(F_{-1}(\bar x^{\J}),\bar z^{\J}))=f_{\circ}(\bar x^{\J},f_{\circ}(f_{-1}(\bar x)^{\J},\bar z^{\J}))=\bar z^{\J}.$\qed
\medskip 
\par Since both $V_{0}^{''}, V_{1}^{''}$ are $\cL$-definable and almost equal to $G^{**}$ (see Claims \ref{$V_{0}^{''}$}, \ref{$V_{1}^{''}$}), $V_{0}^{''}\cap V_{1}^{''}$ is $\cL$-definable and almost equal to $G^{**}$. So it can be expressed as a finite union of continuous correspondences by Proposition \ref{prop:cell}. Let $V_{2}$ be the union of those with the same $\cL$-dimension as $G^{**}$.
\cl \label{$V$}Let $V:=V_{2}\cap F_{-1}(V_{2})$. Then $V$ is $\cL$-definable, relatively open in $G^{**}$  and almost equal to $G^{**}$.
Moreover $F_{-1}$ maps $V$ to $V$ and is continuous.
\ecl
\prcl It remains to show that $F_{-1}$ maps $V$ to $V$. An element of $V$ can be written as $\bar a\in V_{2}$ and $F_{-1}(\bar b)$ for $\bar b\in V_{2}$. Since $V_{2}\subset V_{1}$, we get that $F_{-1}(F_{-1}(\bar b))=\bar b$. 
\qed
\medskip
\par By construction, $V_{2}$ is a relatively open definable large subset of $G^{**}$. We define 
$Y:=\{(\bar x,\bar y)\in (V\times V)\cap Y_{0}:\;F_{\circ}(\bar x,\bar y)\in V\}$.
Note that $Y$ is large in $V\times V$. If we choose $\bar x^{\J}\in V$ generic, then $\bar y^{\J}\in V$ generic over $\bar x^{\J}$.
Then $F_{\circ}(\bar x^{\J},\bar y^{\J})$ is generic over $\bar x^{\J}$ by Claim \ref{ageneric} and so belongs to $V$. Note that on $Y$, the map $F_{\circ}$ is continuous.

\

By Claim \ref{$V$}, $V$ is a large definable relatively open subset of $G^{**}$. By the paragraph above $Y$ is a large definable relatively open subset of $G^{**}\times G^{**}$ and $F_{\circ}$ is continuous on $Y$ (statement (3) of the proposition. Statement (1) of the proposition follows by construction, (2) is Claim \ref{$V$} and 
(4) is proven in the same way as we did with $V_0$ and $Y_0$ (see Claim \ref{$V_{0}$}).
\qed

\medskip

Now let us recall the notion of infinitesimal neighbourhoods in the context of $\cL$-open theories $T$ of topological fields. 

\defn \cite[Definition 2.3]{Pet}\label{infinit} Let $\cM\models T$ and let $a\in M^n$. Then the $M$-infinitesimal neighbourhood of $a\in M^n$ is the partial type consisting of all formulas with parameters in $M$ that define an open subset of $M^n$ containing $a$. Given $\cN$ an $\vert M\vert^+$-saturated extension of $\cM$, we denote by $\mu_{a}(\cN)$ its realization in $\cN$. Let $X$ be a $M$-definable subset of $M^n$ containing $a$, then $\mu_{a}(X)=\mu_{a}(\cN)\cap X(\cN)$. 
\par We will also use the notation $u\sim_M 0$ to mean that $u\in \mu_0(\cN)$ and use the term for such elements $u$, $M$-infinitesimals.
\edefn
If $a$ is generic in $X$, then $\mu_{a}(X)$ is independent of the choice of $X$, namely we have the analog of \cite[Fact 2.4]{Pet}, using property \ref{dim} of the dimension function in models of $T$. For convenience of the reader, we prove it below.
\lmm \cite[Fact 2.4]{Pet} Let $X$ be an $A$-definable subset of $M^n$ and $a$ a generic element of $X$ over $A$. Then for any $A$-definable set $Y$, if $a\in Y$, then $\mu_{a}(X)\subset \mu_{a}(Y)$. In particular if $\dim(X)=\dim(Y)$, then $\mu_{a}(X)=\mu_{a}(Y)$.
\elmm
\pr Consider $X\cap Y$. Then since $a$ is generic in $X$, $a\in Int_{X}(X\cap Y)$. Since the topology is definable, there exists an open $A$-definable set $U$ containing $a$ such that $U\cap X\subset X\cap Y$. So, $U\cap X=U\cap X\cap Y$. It follows that $\mu_{a}(X)=\mu_{a}(X\cap Y)$ and so $\mu_{a}(X)\subset \mu_{a}(Y)$. If $\dim(X)=\dim(Y)$, then $a$ is also generic in $Y$ and the reverse inclusion holds, namely $\mu_{a}(Y)\subseteq \mu_{a}(X)$.
\qed
\medskip
\par The above lemma allows us to introduce the following notation.
\nota \cite[Definition 2.5]{Pet} \label{mu} Given $a\in M^n$, $A\subset M$ and $X$ an $A$-definable subset of $M^n$ with $a$ generic in $X$ over $A$, we denote $\mu^{\cM}(a/A)$ (or $\mu(a/A)$) the set $\mu_{a}(X)$.
\enota

\nota \cite[Definition 2.10]{Pet} Let $\cM\models T$ and let $O$ be an open subset of $M^n$. Let $p_{1}, p_{2}: O\to M$ be two maps and let $y\in O$. Then $p_{1}\sim_{y} p_{2}$
means that for some open neighbourhood $U\subset O$ of $y$, $p_{1}\restriction U=p_{2}\restriction U$.
\enota

\lmm \label{faithful} \cite[Lemma 2.11]{Pet} Let $\cM$ be a sufficiently saturated model of $T$. Let $V\subset M^n$ be a $\cL$-definable open subset in $\cM$; let $\cF:=\{p_b:V\to V\colon  b\in V\}$ be a family of $\cL$-definable bijections of $V$. Let $x$ be a generic element of $V$ and let $a_1\in V$ generic over $\{ x\}$. Then there exist definable open subsets $W$ containing $x$ and $U$ containing $a_1$ such that for every $a_1', a_1''\in U$ and $y\in W$, if $p_{a_1'}\sim_{y} p_{a_{1}''}$, then $p_{a_{1}'}\restriction W=p_{a_{1}''}\restriction W$.
\elmm
\pr Denote by $[a_{1}]_{x}:=\{u\in V\colon p_{u}\sim_{x} p_{a_{1}}\&\; p_{u}\in \cF\}$. Let $a\in [a_{1}]_{x}$ be generic over $\{a_{1}, x\}$. Let $W_{0}$ be an open neighbourhood of $x$ such that $p_{a_{1}}\restriction W_{0}=p_{a}\restriction W_{0}$. Recall that the topology on $K$ is definable and let $\bar d_{0}$ be parameters such that $W_{0}=\bar \chi(M,\bar d_{0})$.
We may assume by shrinking $W_{0}$ if necessary that $\bar d_{0}$ is independent from $a_{1},a, x$. We can express by an $\cL$-formula in $a_{1},a,x, \bar d_{0}$ the following property:
\[p_{a_{1}}\sim_{x} p_{a}\rightarrow p_{a_{1}}\restriction W_{0}=p_{a}\restriction W_{0}.
\] 
So it continues to hold for all $a'$ in an open neighbourhood $U_{0}$ of $a_{1}$ and this last property can be expressed by an $\cL$-formula in $a_{1}$, $x$ and $\bar d_0$:
\[\forall a'\in U_{0}\;(p_{a_{1}}\sim_{x} p_{a'}\rightarrow p_{a_{1}}\restriction W_{0}=p_{a'}\restriction W_{0}).
\]
Since $a_{1}$ and $x$ are independent, this formula continues to hold in an open neighbourhood $U_{1}$ of $a_{1}$ and an open neighbourhood $W_{1}$ of $x$. Then $U:=U_{1}\cap U_{0}$ and $W:=W_{0}\cap W_{1}$ are the sought neighbourhoods.
\qed

\thm \label{local} Let $\cS$ be sorts in  $\cL^{\mathrm{eq}}$ such that $T$ admits elimination of imaginaries in $\cL^{\cS}$. Let $\cK\models T_{\delta}^*$ and assume that $\cK$ is sufficiently saturated. Let $\cG:=(G,f_{\circ} ,f_{-1},e )$ be an $\cL_{\delta}$-definable group in (the field sort of) $\cK$ (over a subset of parameters). Let $G^{**}$ be an $\cL$-open envelop of $G$.
Then there exists a type $\cL$-definable topological group $H$ (over some parameters) with $\dim(H)=\dim(G^{**})$. 
\ethm
\pr We keep the same notations as in Proposition \ref{largeV}; in particular we will use $F_{-1}$, $F_{\circ}$, $V$ and $Y$. In order to simplify notations we will not indicate the parameters we are working with and simply use $\cL$; however
we will indicate all the additional parameters. 
\par Let $a_{1}\in G^{**}$ be $\cL$-generic and let $a_{2}$ be $\cL$-generic over $a_{1}$.
Then by Proposition \ref{largeV} $(4)$, we have: $(a_{2},a_{1})\in Y$ and $(F_{-1}(a_{2}),F_{\circ}(a_{2},a_{1}))\in Y$. By Proposition \ref{largeV} $(3)$, $a_3\in V$.
\cl Let us show that $a_{3}=F_{\circ}(a_{2},a_{1})$ is $\cL$-generic in $V$.
\ecl
\prcl This is similar to the proof of Claim \ref{ageneric}.
We have that 

\noindent $acl(a_2, F_{\circ}(a_2,a_1))\subset acl(F_{\circ}(a_2,a_1), a_2, F_{-1}(a_2))\subset acl(a_2,a_1)\subset acl(F_{\circ}(a_2,a_1),a_2)$. So, by equation (\ref{dim})
\begin{align*}
 \dim(a_1,a_2)=&\dim(a_2/a_1)+\dim(a_1)=\\
 \dim(a_2, F_{\circ}(a_2,a_1))=&\dim( F_{\circ}(a_2,a_1)/a_2)+\dim(a_2).
 \end{align*}
  By assumption $\dim(a_2/a_1)=\dim(G^{**})=\dim(a_2)=\dim(a_1)$.
  
  Therefore,  $F_{\circ}(a_2,a_1)$ is generic over $a_2$ and so generic in $G^{**}$.
\qed

\

So for each $1\leq i\leq 3$, the map $V\to V: u\mapsto F_{\circ}(a_{i},u)$ is a bijection on $V$ by Claim \ref{$V_{1}^{''}$} and the definition of $V$. Note that this is a $\cL$-definable property of each $a_{i}$ and so this will hold for any $a_{i}'\in \mu(a_{i})$, $1\leq i\leq 3$. 

We define a family of bijections on $V$ around each $a_{i}$, $1\leq i\leq 3$, using the following notations.
\begin{enumerate}
\item Denote for $a_{1}'\in \mu(a_{1})$, $p_{a_{1}'}:V\to V: u\mapsto F_{\circ}(a_{1}',u)$ and by \\
$\cP:=\{p_{a_{1}'}:\;a_{1}'\in \mu(a_{1})\}$.
\item Denote for $a_{2}'\in \mu(a_{2})$, $q_{a_{2}'}:V\to V:u\mapsto F_{\circ}(a_{2}',u)$ and by 
\\$\cQ:=\{q_{a_{2}'}:\;a_{2}'\in \mu(a_{2})\}$.
\item Denote for $a_{3}'\in \mu(a_{3})$, $h_{a_{3}'}:V\to V:u\mapsto F_{\circ}(a_{3}',u)$ and by 
\\$\cH:=\{h_{a_{3}'}:\;a_{3}'\in \mu(a_{3})\}$.
\end{enumerate}
Let $x\in V$ be $\cL$-generic over $\{a_{1},a_{2}\}$. In particular $\mu(x)=\mu(x/\{a_{1},a_{2}\})$.
As in the group configuration theorem \cite[Theorem 3.4]{Pet}, we consider the following subgroup $H$ of the group $Sym(\mu(x))$ of permutations on $\mu(x)$ generated by: $\{p_{a_{1}'}^{-1}\circ p_{a_{1}''}:\; a_{1}', a_{1}''\in \mu(a_{1})\}$. 

\cl $\forall p\in \cP\,\forall q\in\cQ\,\exists h\in \cH\,\; q\circ p=h$ on $\mu(x)$.
\ecl 
\prcl Let $a_{1}'\in \mu(a_{1})$ be such that  $p=p_{a_{1}'}$, let $a_{2}'\in \mu(a_{2})$ be such that $q=q_{a_{2}'}$. 
The map $F_{\circ}$ is continuous on $Y$, 
so for every neighbourhood $W$ of $a_{3}$ there exists $W_{1}$ a neighbourhood of $a_{1}$ and $W_{2}$ a neighbourhood of $a_{2}$ such that $F_{\circ}(W_{2},W_{1})\subset W$.
\par First let us show that $F_{\circ}(a_{2},F_{\circ}(a_{1},x))=F_{\circ}(F_{\circ}(a_{2},a_{1}),x)$. By Claim \ref{$V_{0}^{''}$}, we know that $\tilde Y_{0}(x)$ is almost equal to $Y_{0}$. Since $x$ is $\cL$-generic over $\{a_{2},a_{1}\}$, we get that $(a_{2},a_{1})\in \tilde Y_{0}(x)$. So, $(a_{2},F_{\circ}(a_{1},x))\in Y_{0}$ and $F_{\circ}(a_{2},F_{\circ}(a_{1},x))=F_{\circ}(F_{\circ}(a_{2},a_{1}),x)$.
Furthermore, $(a_{2},a_{1})\in \Int_{Y_{0}}(\tilde Y_{0}(x))$. So for $a_{1}'\in \mu(a_{1})$ and $a_{2}'\in \mu(a_{2})$, we get that $(a_{1}',a_{2}')\in Y_{0}(x)$ and so $F_{\circ}(a_{2}',F_{\circ}(a_{1}',x))=F_{\circ}(F_{\circ}(a_{2}',a_{1}'),x)$.
\par Now we can express by a $\cL(\{a_{1},a_{2}\})$-formula (in $x$) the property that $(a_{2},a_{1})\in \Int_{Y_{0}}(\tilde Y_{0}(x))$. So for any   $x'\in \mu(x)=\mu(x/\{a_{1},a_{2}\})$, since $x$ is generic in $V$, we get that $(a_{2},a_{1})\in \Int_{Y_{0}}(\tilde Y_{0}(x'))$.
Therefore, $F_{\circ}(a_{2},F_{\circ}(a_{1},x'))=F_{\circ}(F_{\circ}(a_{2},a_{1}),x')$. By the same reasoning as above, it also holds for $a_{1}', a_{2}'$ in place of $a_{1}, a_{2}$.
\qed

 \cl $\forall h\in \cH\,\forall q\in \cQ\,\exists p\in \cP\; q\circ p=h$ on $\mu(x)$.
 \ecl
 \prcl  Let $a_{2}'\in \mu(a_{2})$ be such that  $q=q_{a_{2}'}$, let $a_{3}'\in \mu(a_{3})$ be such that $h=h_{a_{3}'}$. Then define a map  $p$ on $u\in \mu(x)$ as follows $p(u):=F_{\circ}(F_{\circ}(F_{-1}(a_{2}'),a_{3}'),u)$. Note that $F_{-1}(V)=V$ and $F_{\circ}(a_{2},F_{\circ}(F_{-1}(a_{2}),a_{3}))=a_{3}$ by Claim \ref{$V_{1}^{''}$} ($a_{2}$ is $\cL$-generic).
Since $F_{\circ}$ and $F_{-1}$ are continuous, it holds in a neighbourhood of $a_{2}$, respectively $a_{3}$. 

 \qed
 
 \cl $\forall h\in \cH\,\forall p\in \cP\,\exists q\in \cQ\; q\circ p=h$ on $\mu(x)$.
 \ecl
 \prcl Let $a_{1}'\in \mu(a_{1})$ be such that  $p=p_{a_{1}'}$, let $a_{3}'\in \mu(a_{3})$ be such that $h=h_{a_{3}'}$. Then define a map $q$ on $u\in \mu(F_{\circ}(a_{1},x)$ by $q(u):=F_{\circ}(F_{\circ}(a_{3}',F_{-1}(a_{1}')),u)$. Note that $F_{-1}(V)=V$ and $F_{\circ}(F_{\circ}(a_{3},F_{-1}(a_{1})),a_{1})=a_{3}$ by Claim \ref{$V_{1}^{''}$} ($a_{1}$ is $\cL$-generic and $a_{3}$ is $\cL$-generic over $a_{1}$).
Since $F_{\circ}$ and $F_{-1}$ are continuous, it holds in a neighbourhood of $a_{1}$, respectively $a_{3}$. 
 \qed
 \cl $\forall p_{1}\in \cP\,\forall p_{2}\in \cP\,\forall p_{3}\in \cP\,\exists p_{4}\in \cP\,\;p_{1}^{-1}p_{2}=p_{3}^{-1}p_{4}$ on $\mu(x)$.
 \ecl
 \prcl We apply the previous claims in order to write first $p_{3}$ as $q_{3}^{-1}\circ h_{3}$ with $q_{3}\in \cQ$ and $h_{3}\in \cH$ by Claim 3.8.
 Then write $p_{1}$ as $q_{1}^{-1}\circ h_{3}$ with $q_{1}\in \cQ$ by Claim 3.10 and finally $p_{2}=q_{1}^{-1}\circ h_{2}$ with $h_{2}\in \cH$ by Claim 3.8.
 Composing these maps, we get $p_{3}p_{1}^{-1}p_{2}=q_{3}^{-1}h_{2}\in \cP$ by Claim 3.9.
   \qed
   
 \cl \label{can} Let $a$ be a fixed element of $\mu(a_{1})$ $\cL$-generic in $V$ over $a_1$. Then group $H$ is equal to $\{p_{a}^{-1}\circ p_{a_{1}'}:\; a_{1}'\in \mu(a_{1})\}$.
  \ecl
 \prcl We apply the previous claim with $p_{3}=p_{a}$.
 \qed 

\cl \label{trans} The action of the group $H$ is transitive on $\mu(x)$.
\ecl
\prcl Since $x\in V$, $F_{-1}(x)\in V$ and since $x$ is generic over $a_{1}$, $(a_{1},x)\in Y$. By Claim \ref{$V_{0}^{''}$}, $F_{\circ}(F_{\circ}(a_{1},x),F_{-1}(x))=a_{1}$. 
So given $x_{2}'\in \mu(x)$ and $x_{3}'\in \mu(F_{\circ}(a_{1},x))$, we may define $a_{1}':=F_{\circ}(x_{3}',F_{-1}(x_{2}'))$ and this element belongs to $\mu(a_{1})$. 

Then we re-apply the same reasoning to $x_{2}''\in \mu(x)$ and $x_{3}'$, and define $a_{1}'':=F_{\circ}(x_{3}',F_{-1}(x_{2}''))$. Then $a_{1}''\in \mu(a_{1})$.
Finally consider $p_{a_{1}''}^{-1}\circ p_{a_{1}'}$; by Claim \ref{can}, it belongs to $H$ and $p_{a_{1}''}^{-1}\circ p_{a_{1}'}(x_{2}')=x_{2}''$.
\qed

\

Finally by Lemma \ref{faithful}, there is a $\cL$-definable open neighbourhoods $W\subset V$ of $x$ and $U\subset V$ of $a_{1}$ such that if $p_{a_{1}'}\sim_{x} p_{a_{1}''}$, then $p_{a_{1}'}\restriction W=p_{a_{1}''}\restriction W$. For $a_{1}'\in U$, let $[a_{1}']_{x}:=\{a_{1}^{''}\in U\colon p_{a_{1}'}\sim_{x} p_{a_{1}''}\}= \{a_{1}^{''}\in U\colon p_{a_{1}'}\restriction W= p_{a_{1}''}\restriction W\}$. Now suppose that $p_{a}^{-1}p_{a_{1}}'(u)=p_{a}^{-1}p_{a_{1}}'(u)$, namely $F_{\circ}(F_{-1}(a),F_{\circ}(a_{1}',u))=F_{\circ}(F_{\circ}(F_{-1}(a),a_{1}''),u)$.
So $F_{circ}(a,F_{\circ}(F_{-1}(a),F_{\circ}(a_{1}',u)))=F_{\circ}(a,F_{\circ}(F_{-1}(a),F_{\circ}(a_{1}',u)))$ which implies that $F_{\circ}(a_{1}',u)=F_{\circ}(a_{1}',u)$. 
So, two elements of $H$ are equal on a neighbourhood of $x$, then they coincide on $U$ by Lemma \ref{faithful}. This will allow us to  identify $H:=\{p_{a}^{-1}\circ p_{a_{1}'}:\; a_{1}'\in \mu(a_{1})\}$ with $\mu(a_1)/E$ where $\mu(a_1)$ is type-definable over $M$ and $E$ is a definable equivalence relation on $V$, defined by $E(a_{1}',a_{1}'')$ if and only if $p_{a_{1}'}\restriction W=p_{a_{1}''}\restriction W$. 
\par The group law is definable since the action of the elements of $H$ on $\mu(x)$ is defined using the $\cL$-definable functions $F_{\circ}$ and $F_{-1}$.
Finally since $H$ acts transitively on $\mu(x)$ (Claim \ref{trans}) its dimension is equal to $\dim(x)=\dim(V)=\dim(G^{**})$.
Furthermore $H$ a topological group, namely the group laws are continuous. (It follows from the fact that $F^{-1}$ is continuous on $V$ and $F_{\circ}$ is continuous on $Y$.)
\qed.

\rem Note that we may associate with $H$ a $\dim(G^{**})$-group configuration \cite[Definition 3.1]{Pet}.
\erem

\section{Annex: Largeness}\label{annex}

\par In the following proposition, we put an extra-hypothesis besides largeness on a differential topological field that we will call c-largeness, in order to embed a differential topological field $(K,\delta)$ endowed with a definable topology into a differential extension model of the scheme (DL). 

A topological field $K$ is c-large if it is large and the following holds. We consider an embedding $K\hookrightarrow K((s_{0},\cdots,s_{n-1}))\hookrightarrow K^*$ with $K^*$ a non principal ultrapower of $K$, $s_0,\ldots,s_{n-1}$ algebraically independent over $M$ and $s_i\sim_K 0$, $0\leq i\leq n-1$ (see Definition \ref{infinit}). We endow $K((s_0,\ldots,s_{n-1}))$ with a valuation $v$, trivial on $K$ and with $0<v(s_{0})<v(s_{1})<\cdots<v(s_{n-1})$. Let $\cM_v$ be the maximal ideal of $K[[s_0,\ldots,s_{n-1}]]$, consisting of elements of strictly positive valuation. Then the additional condition we ask is the following: $\cM_v$ embeds into the $K$-infinitesimal elements of $K^*$. 
\begin{proposition} Let $\cK_\delta$ be a c-large topological differential field. 
Then we can embed $\cK_\delta$ into a model of the scheme $(DL)$ and with the same $\cL$-theory as $\cK$.
\end{proposition}
\pr We only show the main step.  As classically done in the construction of existentially closed models, one first enumerate the existential formulas (with parameters in the ground field $K$, belonging to a certain family) one wants to satisfy in an extension and then one redo the construction $\omega$ times in order to get an existentially closed model (with respect to that family of formulas) containing $K$ as the union of elementary extensions of $K$. 
\par  Consider a differential polynomial $p(x)$ in one variable in $K\{x\}$ of order $n>0$ such that $p^*(\bar a)=0\,\&\,s_{p}^*(\bar a)\neq 0$ for some $\bar a:=(a_0,\ldots,a_n)\in K^{n+1}$. Then we want to find an $\cL$-elementary extension of $\cK$ and an extension of $\delta$ such that in that elementary extension there is a differential solution of $p(x)=0$ with $\bar \delta^n(b)\sim_K \bar a$. 
We first consider an ultrapower $\cK^*$ of $\cK$ which is $\vert K\vert^+$-saturated and in that ultrapower $n$ elements $s_0,\ldots,s_{n-1}$ with $s_0\sim_K 0$ and for $i\geq 1$, $s_i\sim_{K(s_0,\ldots,s_{i-1})} 0$. Then $s_0,\ldots, s_{n-1}$ are algebraically independent over $K$. Since $\cK$ is large, we have that $\cK\subset \cK((s_0,\ldots,s_{n-1}))\subset \cK^{**}$
for some ultrapower of $\cK$. Furthermore, we may assume that $s_i$, $0\leq i\leq n-1$, are still $K$-infinitesimals in $K^{**}$.
\par Set $\bar s_{n}:=(s_{0},\cdots, s_{n-1})$. 
We rewrite  $\bar a$ as $(\bar a_n, a_n)$ with $\vert \bar a_n\vert =n$, $\vert a_n\vert=1$ and $p^*(\bar x)$ as $p^*(\bar x_n, x_n)$ with $\vert \bar x_n\vert=n$, $\vert x_n\vert=1$. We consider the polynomial in $x_{n}$: $p^*(\bar a_n+\bar s_{n}, x_n)$. We endow the field $K((s_{0},\ldots,s_{n-1}))$ with a valuation $v$, trivial on $K$ and with $0<v(s_{0})<v(s_{1})<\cdots<v(s_{n-1})$. (So the value group is isomorphic to $\Z^n$ with the lexicographic order.) We are in the conditions of Hensel's Lemma since $v(p^*(\bar a_{n}+\bar s_{n},a_{n}))>0$ and $v(s_{p}^*(\bar a_{n}+\bar s_{n},a_{n}))=0$. So there is $b\in K((s_{0},\cdots,s_{n-1}))$ such that $v(b-a_n)>0$ and 
\begin{equation}\label{1}
p^*(\bar a_{n}+\bar s_{n},b)=0\,\&\,s_{p}^*(\bar a_{n}+\bar s_{n},b)\neq 0.
\end{equation}
By c-largeness, we have $b\sim_{K} a_{n}$. 
We extend $\delta$ on $K(s_{0},\cdots,s_{n-1})$ by setting $\delta_{1}(a_i+s_{i})=a_{i+1}+s_{i+1}$, $0\leq i<n-1$ and $\delta_{1}(a_{n-1}+s_{n-1})=b\in K^*$.
Thanks to equation (\ref{1}), $\delta_{1}(b)$ is completely determined and belongs to the subfield of $K(s_{0},\cdots,s_{n-1})$ of $K^*$.
We set $K_{1}=K(s_{0},\cdots,s_{n-1})$ and extend $\delta_1$ on the relative algebraic closure of $K_1$ in $K^*$. We apply Lowenheim-Skolem theorem to embed $K_1$ in an $\cL$-elementary extension $\tilde \cK$ of $\cK$ of the same cardinality and extend $\delta_1$ on $\tilde K$.
\qed

\end{document}